\documentclass[review]{elsarticle}

\usepackage{lineno,hyperref}
\usepackage{amsmath}
\usepackage{amssymb}
\usepackage{booktabs}
\usepackage{mathtools}
\usepackage{caption}
\usepackage{subcaption}
\usepackage{tabularray}
\usepackage{graphicx}
\usepackage[export]{adjustbox}
\usepackage{booktabs}
\usepackage{chngcntr}
\usepackage{algorithm,algorithmic}
\usepackage{xcolor}
\usepackage{setspace}
\usepackage{multirow}
\modulolinenumbers[255]

\usepackage[inkscapelatex=false]{svg}
\newcommand{\phiDef}{\boldsymbol{\phi}_{\mathbf{d}}}
\newcommand{\phiRB}{\boldsymbol{\phi}_{\mathbf{r}}}
\newcommand{\psiDef}{\boldsymbol{\psi}_{\mathbf{d}}}
\newcommand{\psiRB}{\boldsymbol{\psi}_{\mathbf{r}}}
\newcommand{\polDeg}{t}

\usepackage{geometry}
\journal{Computer Methods in Applied Mechanics and Engineering}

\begin{document}

\begin{frontmatter}

\title{Accelerating structural optimization via EIFEM: a ROM-based preconditioner}

\author[CIMNE,RMEE]{R.Rubio\corref{mycorrespondingauthor}}
\author[CIMNE,ESEIAAT]{A.Ferrer}
\author[CIMNE,RMEE]{J.A. Hernández}
\author[EPFL]{P.Antolin}

\cortext[mycorrespondingauthor]{ Corresponding author \\ E-mail address: raul.rubio.serrano@upc.edu (R.Rubio) \\ }

\address[CIMNE]{Centre Internacional de Mètodes Numèrics en Enginyeria (CIMNE),C/ Gran Capitán
S/N UPC Campus Nord, Edifici C1, 08034 Barcelona, Spain.}
\address[RMEE]{E.S. d’Enginyeries Industrial, Aeroespacial i Audiovisual de Terrassa (ESEIAAT), Departament de Resistència de Materials i Estructures a l’Enginyeria (RMEE), Universitat Politècnica de Catalunya (UPC), C/ Colom 11, 08222 Terrassa, Spain}
\address[ESEIAAT]{E.S. d’Enginyeries Industrial, Aeroespacial i Audiovisual de Terrassa (ESEIAAT), Departament de Física (FIS), Universitat Politècnica de Catalunya (UPC), C/ Colom, 11, 08222 Terrassa,
Spain.}
\address[EPFL]{Institute of Mathematics, École polytechnique fédérale de Lausanne, Station 8, CH-1015 Lausanne, Switzerland}





\begin{abstract}
Topology optimization (TopOpt) is a standard tool for structural conceptual design, providing optimal structures with great design freedom. However, it has a high computational costs due to the repeated evaluation of high-fidelity finite element models and frequently produces complex geometries that require extensive interpretation for manufacturability. As an alternative, this work presents a parameterized optimization framework based on Reduced Order Modeling (ROM) and preconditioning techniques. We employ the Empirical Interscale Finite Element Method (EIFEM) coupled with the Discrete Empirical Interpolation Method (DEIM) to construct localized, parameter-dependent reduced operators. From an optimization perspective, this approach can be interpreted as a preconditioning strategy in which inexact gradients are used to accelerate convergence. Furthermore, parameterizing the design space in terms of explicit geometric features, such as inclusion radii or lattice widths, significantly improves the manufacturability of the resulting optimized designs. We assess the performance of the method in terms of computational time, design topology and structural performance with TopOpt as baseline for three different unit cell geometries and three different benchmarks. While the parametric ROM framework requires an initial offline training phase, our comparative analysis demonstrates that it drastically accelerates the online optimization loop and that it can produce even stiffer design in some of our experiments.

\end{abstract}


\end{frontmatter}

\nolinenumbers

\section{Introduction}

Structural optimization relies on iterative algorithms to distribute material within a design domain, minimizing an objective function subject to physical and/or geometric constraints \cite{BendsoeSigmund2003,Sigmund2013}. Over the past decades, topology optimization (TopOpt) has become the standard framework for conceptual design. By defining a characteristic function that indicates the presence or absence of material, TopOpt seeks the optimal layout without a priori assumptions about the structural topology. Since the binary  problem is ill-posed and non-differentiable \cite{Allaire2007}, it is conventionally relaxed using density-based \cite{BendsoeSigmund2003} methods or level-set functions \cite{Allaire2004}. In the density approach, a continuous variable defines the material distribution, and the optimization problem is solved using gradient-based algorithms such as the Method of Moving Asymptotes (MMA) \cite{Svanberg1987MMA} or the Null space algorithm \cite{Feppon2020}.

While TopOpt provides maximum design freedom, it suffers from severe computational bottlenecks. Evaluating the objective function and its sensitivities requires solving a high-fidelity finite element model at every iteration \cite{Allaire2004}. For high-resolution 3D domains, repeatedly solving large-scale linear systems becomes computationally prohibitive. Furthermore, the resulting designs often feature complex, free-form topologies that require post-processing and interpretation to ensure manufacturability, even with advanced techniques such as additive manufacturing. In order to increase manufacturability, minimum length scale \cite{Lazarov2016,Amstutz2022} and overhang constraints \cite{Allaire2017,Torres2025} are introduced. For an overview of the different difficulties arising in additive manufacturing and strategies to tackle those, the reader is referred to \cite{Bayat2023}.

Parametric optimization offers an alternative by restricting the design space to a finite set of explicit geometric features. Instead of operating on a grid of independent density variables, the physical system is governed by a parameter vector that represents variations of the thickness of a certain component or the radius of internal inclusions, among others. This a priori restriction reduces the design space dimensionality and ensures directly manufacturable outcomes. However, the requirement to assemble and solve the full-order state equation at every iteration remains the primary computational barrier.

To lift the computational burden during the optimization loop, Reduced Order Modeling (ROM) techniques construct low-dimensional surrogates of high-fidelity models. Projection-based ROMs approximate the solution manifold using a reduced basis, typically extracted via Proper Orthogonal Decomposition (POD) from a set of training snapshots \cite{Berkooz1993}. If the underlying partial differential equation (PDE) exhibits affine parameter dependence, the reduced operators can be pre-computed in an offline phase. The online optimization then requires only trivial scalar summations and the inversion of a small, dense system, effectively decoupling the computational cost from the full mesh size \cite{quarteroni2016reduced}.

Despite its efficiency, the standard offline-online ROM strategy fails when the parameterization alters the domain geometry. Geometric variations induce non-affine dependencies in the discretized operators. Consequently, evaluating the weak form integrals requires assembling the full high-fidelity matrices at every iteration, negating the computational advantages of the parameterized model. To recover the offline-online decoupling for non-affine parametric PDEs, hyper-reduction techniques such as the Discrete Empirical Interpolation Method (DEIM) \cite{Chaturantabut2010,quarteroni2016reduced,NEGRI2015} are necessary. DEIM approximates the non-affine terms by interpolating them at a small subset of spatial locations, known as interpolation indices. By enforcing exact interpolation at these specific points, DEIM reconstructs the reduced operators without traversing the full high-fidelity mesh.

Global ROMs often struggle to capture localized parameter variations, as modifying a local geometric feature requires updating the global reduced basis \cite{Buhr2021}. Domain decomposition strategies address this limitation by constructing local reduced operators. The Empirical Interscale Finite Element Method (EIFEM), for example, partitions the global domain into subdomains and computes localized reduced bases. EIFEM introduces an inter-scale operator that maps interface displacements to the subdomain interior, allowing the construction of a coarse stiffness matrix. In a parametric setting, these local inter-scale operators and coarse matrices can be parameterized directly using DEIM as a function of the geometric variables. This modularity allows localized geometric updates to be evaluated independently, avoiding global reassembly and offering a highly scalable framework for structural optimization.

A similar approach can be found in \cite{Diercks2025}, where the authors randomly sample the space to obtain a representation of the subdomain displacements; or in \cite{NEZDYUR2026} where the coarse displacements are applied directly on the domain boundary and harmonically extended. In our experience, as discussed for the auxetic cell in \cite{RUBIO2025}, applying the coarse scale displacements on the boundary of a single domain may artificially over-stiff the ROM, especially for units cells whose edges are discontinuous. Other recent works related to the optimization of lattice structures based on domain decomposition \cite{TosselliWidlund2006,Doolean2015,Farhat2001FETIDP,mandel2008multispace} and/or multiscale approaches can be found in \cite{Joskova_Tyburec_Doskar_2025,hirschler2026,Wang2023,Zhou2026,KIM2024} and the references therein.

This work systematically compares ROM-parameterized optimization against standard topology optimization. We investigate the trade-offs between design space restriction, offline training costs, and online computational efficiency. Specifically, we employ a parameterized EIFEM framework where the inter-scale mapping and local coarse operators are interpolated via DEIM over geometric parameter spaces, such as inclusion radii or lattice width, for architected materials.  This strategy can be regarded as a preconditioner for the optimization loop by means of an inexact gradient \cite{Nocedal2006}, where the solution of the PDE is approximated by EIFEM. We contrast this approach with density-based TopOpt, evaluating both methods on their convergence behavior, the computational time per iteration, and the structural performance of the optimized designs. By benchmarking the finite-dimensional parametric ROM against the topology optimization framework, we establish guidelines for selecting the appropriate optimization strategy based on available computational budgets and manufacturing constraints.

The remainder of this article is organized as follows. Section \ref{sec:EIFEM} briefly describes the EIFEM, starting with the Localized Lagrange Multiplier (LLM) method \cite{Park2000,Brezzi2005} and then introducing the dimensionality reduction used to obtain the coarse stiffness matrix of a given subdomain. This section also describes the preconditioning strategy developed in \cite{RUBIO2025}. Section \ref{Section: Parameterized EIFEM} describes and validates the parameterization of the interscale and coarse stiffness operators using DEIM, followed by high-order Lagrange polynomial interpolation to ensure a continuous dependence of the EIFEM operators on the geometrical parameters required by the optimization loop. Section \ref{Section:optimization} then presents the formulation of both density-based topology optimization and EIFEM-based optimization. The numerical assessment of the proposed methodology is presented in Section \ref{Results}. Finally, Section \ref{Conclusions} summarizes the main conclusions and discusses directions for future research.


\section{EIFEM: Empirical Interscale Finite Element Method} \label{sec:EIFEM}

The EIFEM leverages domain decomposition strategies to construct local reduced operators. To illustrate the method, let us briefly recap the formulation introduced in \cite{RUBIO2025}. Let us describe the variational formulation of a partitioned body using the LLM method \cite{Park2000,Brezzi2005}. Consider a domain $\Omega \subset \mathbb{R}^d$ (d=2 or 3 for 2D or 3D, respectively) divided into two non-overlapping subdomains $\Omega_1 \cap \Omega_2 = \emptyset$ such that $\bar{\Omega} = \bar{\Omega}_1 \cup \bar{\Omega}_2$. We denote the boundary of $\Omega$ as $\partial\Omega$ and the boundaries of the subdomains $\Omega^i$ as $\partial\Omega^i$. The set $\Gamma := \bigcup_{i} (\partial\Omega^i \setminus \partial\Omega)$, refers to the internal skeleton. For simplicity, the material is assumed to be in the linear elastic regime. We could then cast the problem as the minimization of the total potential energy subjected to the continuity constraint:

\begin{equation} \label{energy functional LLM}
 \begin{aligned}
\underset{u_i}{\min} \ &  \sum_{i=1}^{2}\frac{1}{2}\int _{\Omega^i} \nabla ^{s}u_i :\mathbb{C}: \nabla ^{s}u_i \ d\Omega \ - \int _{\Omega^i}u_i \cdot  b_i \ d\Omega \ -
\int _{\partial\Omega^i \setminus \Gamma}u{_i} \cdot t_i \ d\Gamma \\
 & \text{s.t.}\ \sum_{i=1}^2\int _{\Gamma} \delta {\lambda_i}\cdot (u_i  - u{_{\Gamma }}) \ d\Gamma  = 0 \quad \forall \delta {\lambda_i},
 \end{aligned}
\end{equation}

\noindent where $u_i$ are the subdomain displacements, $u_\Gamma$ the interface displacement, $\lambda$ the Lagrange multipliers, $b_i$ the body forces and $t_i$ the traction forces. We denote by $\mathbb{C}$ the fourth order elasticity tensor and $\nabla^s$ the symmetric gradient operator.

After a suitable finite element discretization, let 
$V_h^i \subset H^1(\Omega^i)$, 
$V_{\Gamma,h} \subset H^{1/2}(\Gamma)$, and 
$M_h^i \subset H^{-1/2}(\Gamma)$ denote the discrete functional spaces and $\mathbf{N}_i: \Omega^i \rightarrow \mathbb{R}^{d \times N_{loc}}$, $\mathbf{N}_\Gamma : \Gamma \rightarrow \mathbb{R}^{d \times n_b}$, and $\mathbf{N}_i^b : \partial \Omega^i \cap\Gamma \rightarrow \mathbb{R}^{d \times n_b}$ be the finite element shape functions spanning the spaces $V_h^i$, $V_{\Gamma,h}$, and $M_h^i$, respectively. Here, $N_{loc}$ is the number of fine-scale DoFs in the subdomain while $n_b$ is the number of fine-scale DoFs associated with the interface (i.e., on $\Gamma$). The corresponding discrete fields are then given by
\begin{equation} \label{3fields approximation}
\setlength{\jot}{2pt}
    \begin{aligned}
  &u_i^h(x)=\mathbf N_i(x)\mathbf u_i,\hspace{14 pt} x\in\Omega^i,\ i=1,2, \\ 
  & u_\Gamma^h(x)=\mathbf N_\Gamma(x)\mathbf u_\Gamma,\hspace{10 pt} x\in\Gamma,  \\
 &\lambda_i^h(x)=\mathbf N_i^b(x)\boldsymbol\lambda_i, \hspace{13 pt} x\in\partial\Omega^i\cap\Gamma,\ i=1,2,
\end{aligned}
\end{equation}

\noindent where the superscript $\mathbf{b}$ indicates restriction to the boundary and $x$ is the position vector.
Inserting these approximations into the energy functional \eqref{energy functional LLM}
and taking variations with respect to the discrete unknowns
$(\mathbf u_1,\mathbf u_2,\boldsymbol\lambda_1,\boldsymbol\lambda_2,\mathbf u_\Gamma)$
yields the discrete stationarity conditions

\begin{equation} \label{System LLM}
\begin{bmatrix}
\mathbf{K_1}       & 0            & \mathbf{B_1A}_{\mathbf{1}}^T  & 0                 & 0\\
0                  & \mathbf{K_2} & 0                  & \mathbf{B_2A}_{\mathbf{2}}^T & 0\\
\mathbf{A_1B}_{\mathbf{1}}^T  & 0            & 0                  & 0                 & -\mathbf{L_1}\\
0                  & \mathbf{A_2B}_{\mathbf{2}}^T & 0                  & 0                 & -\mathbf{L_2}\\
0                  & 0            & -\mathbf{L}_{\mathbf{1}}^T    & -\mathbf{L}_{\mathbf{2}}^T   & 0
\end{bmatrix} \ 
\begin{bmatrix}
\mathbf{u_1}\\
\mathbf{u_2}\\
\boldsymbol{\lambda_1}\\
\boldsymbol{\lambda_2} \\
\mathbf{u_\Gamma}
\end{bmatrix} =
\begin{bmatrix}
\mathbf{F_1}\\
\mathbf{F_2}\\
\mathbf 0\\
\mathbf 0\\
\mathbf 0
\end{bmatrix} ,
\end{equation}

\noindent where we made use of the following definitions:

\begin{equation}\label{operators}
\begin{aligned}
\mathbf{K}_i &= \int_{\Omega^i} \nabla^{s}\mathbf{N}_i : \mathbb{C} : \nabla^{s}\mathbf{N}_i \, d\Omega,
\qquad
\mathbf{A}_i = \int_{\Gamma} \mathbf{N}_i^{b^{T}}\mathbf{N}_i^b \, d\Gamma, \\
\mathbf{F}_i &= \int_{\Omega^i} \mathbf{N}_i^{T} b_i \, d\Omega
+ \int_{\partial\Omega^i \setminus \Gamma} \mathbf{N}_i^{b^{T}} t_i \, d\Gamma,
\qquad
\mathbf{L}_i = \int_{\Gamma} \mathbf{N}_i^{b^{T}}\mathbf{N}_{\Gamma} \, d\Gamma,
\quad (i=1,2).
\end{aligned}
\end{equation}

\noindent and, following \cite{DVORAK2024}, $\mathbf{B_i}$ is a boolean assembly operator of the subdomain $i$. As shown in \cite{RUBIO2025}, given $\mathbf{u_\Gamma}$, the system \eqref{System LLM} is separable 

 \begin{equation} \label{separable 1}
\begin{bmatrix}
\mathbf{K_1}       & \mathbf{B_1A}_{\mathbf{1}}^T \\
\mathbf{A_1B}_{\mathbf{1}}^T  & 0                 
\end{bmatrix} 
\begin{bmatrix}
\mathbf{u_1}\\
\boldsymbol{\lambda_1}\\
\end{bmatrix} =
\begin{bmatrix}
\mathbf{F_1}\\
\mathbf{L_1}\mathbf{u_\Gamma}
\end{bmatrix}
\quad ; \quad 
\begin{bmatrix}
\mathbf{K_2}       & \mathbf{B_2A}_{\mathbf{2}}^T \\
\mathbf{A_2B}_{\mathbf{2}}^T   & 0                 
\end{bmatrix} 
\begin{bmatrix}
\mathbf{u_2}\\
\boldsymbol{\lambda_2}\\
\end{bmatrix} =
\begin{bmatrix}
\mathbf{F_2}\\
\mathbf{L_2}\mathbf{u_\Gamma}
\end{bmatrix}
\end{equation}

\noindent and allows us to construct reduced local operators for every subdomain. Hence in the rest of this section we \emph{drop the indices} and focus only on one subdomain. 

At this point, we introduce the affine decomposition approximation
\begin{equation} \label{reduced_basis}
\setlength{\jot}{2pt}
\begin{aligned}
  &{u^h(x)} \approx \phiRB(x)\mathbf{u_r} +\phiDef (x) \mathbf{u_d} , \hspace{16pt}  x \in \Omega  , \\ 
  &\lambda^h(x) \approx \psiRB (x) \boldsymbol{\lambda_r} + \psiDef (x) \boldsymbol{\lambda_d} , \hspace{13pt} x \in \partial\Omega \cap\Gamma  ,  \\
 &{u^h_{\Gamma}(x)} \approx \boldsymbol{\mathit{V}}(x) \mathbf{u_\Gamma} , \hspace{65pt} x \in \partial\Omega \cap\Gamma ,
\end{aligned}
\end{equation}

\noindent where the subscripts $\mathbf r$ and $\mathbf d$ denote rigid body and deformational components, respectively, and $\boldsymbol{\mathit{V}}$ stands as a priori chosen set of coarse functions approximating the kinematics of the interface. In this work, we choose $\boldsymbol{\mathit{V}}$ to be coarse FE functions. Each basis in the decomposition is then discretized using the same FE basis as in Equation \eqref{3fields approximation}

\begin{equation} \label{FE discretization}
\setlength{\jot}{2pt}
\begin{gathered}
\boldsymbol{\phi_d}(x)    = \mathbf N(x) \mathbf{\Phi_{d}} \quad ; \quad \boldsymbol{\phi_r}(x)  = \mathbf N(x) \mathbf{\Phi_{r}}  ,
\\
\boldsymbol{\psi_d}(x)   = \mathbf{N^{b}}(x) \mathbf{\Psi_{d}}  \quad ; \quad \boldsymbol{\psi_r}(x) = \mathbf{N^{b}}(x) \mathbf{\Psi _{r}}  , 
\\
\boldsymbol{\mathit{V}}(x) \ = \mathbf{N_{\Gamma }}(x) \mathbf{V_{\Gamma }}
\end{gathered}
\end{equation}

\noindent where $\mathbf{\Phi _{d}} \in \mathbb{R}^{N_{loc} \times p}$ ($p \ll N_{loc}$), $\mathbf{\Phi _{r}}  \in \mathbb{R}^{N_{loc} \times n_{\text{rb}}}$ (with $n_{\text{rb}} = 3$ in 2D and $n_{\text{rb}} = 6$ in 3D), $\mathbf{\Psi _{d}} \in \mathbb{R}^{n_b \times p}$, $\mathbf{\Psi _{r}}\in \mathbb{R}^{n_b \times n_{\text{rb}}}$, and $\mathbf{V_{\Gamma}}\in \mathbb{R}^{n_b \times m}$ being $p$ the number of deformational modes, $n_{\text{rb}}$ the number of rigid body modes and $m$ the chosen number of interface modes. 


Plugging in Equations (\ref{reduced_basis}) and (\ref{FE discretization}) into Equation (\ref{separable 1}) we get the reduced version of the LLM system:
\begin{equation}\label{reduced system}
\begin{bmatrix}
\mathbf{\mathbf{K_{dd}}} & 0 & \mathbf{A_{dr}}  & \mathbf{A_{dd}} & 0 \\
0                        & 0 & \mathbf{A_{rr}}  & 0                & 0 \\
\mathbf{A}_{\mathbf{dr}}^T  & \mathbf{A}_{\mathbf{rr}}^T  & 0 & 0 &  -\mathbf{L_{rv}} \\
\mathbf{A}_{\mathbf{dd}}^T       & 0 & 0                 & 0                & -\mathbf{L_{dv}}  \\
0                        & 0 & -\mathbf{L}_{\mathbf{rv}}^T & -\mathbf{L}_{\mathbf{dv}}^T & 0
\end{bmatrix}  
\begin{bmatrix}
\mathbf{u_d}\\
\mathbf{u_r}\\
\boldsymbol{\lambda_r}\\
\boldsymbol{\lambda_d} \\
\mathbf{u_\Gamma}
\end{bmatrix} =
\begin{bmatrix}
\mathbf{\Phi}_{\mathbf{d}}^T \mathbf{F}\\
\mathbf{\Phi}_{\mathbf{r}}^T \mathbf{F}\\
0\\
0\\
0
\end{bmatrix} 
\end{equation}

\noindent where we used the following definitions
\begin{equation}\label{reduced_operators}
\setlength{\jot}{2pt}
\begin{aligned}
\mathbf{K_{dd}} &:= \mathbf{\Phi}_\mathbf{d}^{T}\,\mathbf{K}\,\mathbf{\Phi_d}, \\[2pt]
\mathbf{A_{dd}} &:= \mathbf{\Phi}_\mathbf{d}^{T}\,\mathbf{BA}^T\,\mathbf{\Psi_d}, 
&\quad
\mathbf{A_{dr}} &:= \mathbf{\Phi}_\mathbf{d}^{T}\,\mathbf{BA}^T\,\mathbf{\Psi_r}, 
&\quad
\mathbf{A_{rr}} &:= \mathbf{\Phi}_\mathbf{r}^{T}\,\mathbf{BA}^T\,\mathbf{\Psi_r}, \\[2pt]
\mathbf{L_{dv}} &:= \mathbf{\Psi}_\mathbf{d}^{T}\,\mathbf{L}\,\mathbf{V_\Gamma}, 
&\quad
\mathbf{L_{rv}} &:= \mathbf{\Psi}_\mathbf{r}^{T}\,\mathbf{L}\,\mathbf{V_\Gamma},
\end{aligned}
\end{equation}

\noindent and the fact that $\mathbf{\Phi}_\mathbf{r}$ lies in the kernel of $\mathbf{K}$, i.e. $\mathbf{K}\mathbf{\Phi}_\mathbf{r} =\mathbf{\Phi}_\mathbf{r}^{T} \mathbf{K} =\mathbf{0}$ .

In Appendix A of \cite{RUBIO2025} we show how to compute the Schur complement of Equation (\ref{reduced system}) to get the equilibrium at the interface in the reduced space
\begin{equation} \label{coarse problem}     
       \mathbf{K_c} \mathbf{u_\Gamma}   =  \boldsymbol{T}^T  \mathbf{F } 
\end{equation}
\noindent where $\boldsymbol{T:=  T_d + T_r}  \in \mathbb{R}^{N_{loc} \times m}$ is an inter-scale operator formed by a deformational component $\boldsymbol{T_d}$ and a rigid body counterpart $\boldsymbol{T_r}$
\begin{equation} \label{operatorT}
\setlength{\jot}{2pt}
\begin{aligned}
 \boldsymbol{T_d} & := 
 \mathbf{\Phi_d} (\mathbf{A}_{\mathbf{dd}}^{T})^{-1}  \mathbf{L_{dv}} =\mathbf{\Phi_d}(\mathbf{\Phi}_\mathbf{d}^{T}  \mathbf B^T \mathbf{\Phi}_\mathbf{d})^{-1}(\mathbf{\Psi}_\mathbf{d}^{T}\mathbf{LV}) \\  
 \boldsymbol{T_r} & := \mathbf{\Phi_r}(\mathbf{A}_{\mathbf{rr}}^{T})^{-1}(\mathbf{L_{rv}}-\mathbf{A}_{\mathbf{dr}}^{T} (\mathbf{A}_{\mathbf{dd}}^{T})^{-1} \mathbf{L_{dv}}) \\
 & = \mathbf{\Phi_r}(\mathbf{\Psi}_\mathbf{r}^{T}
 \mathbf A \mathbf B^T \mathbf{\Phi_r})^{-1}
\mathbf{\Phi}_\mathbf{r}^{T}  \mathbf B \mathbf A ( \mathbf I - \mathbf B^T \mathbf{\Phi_d} (\mathbf{\Psi}_\mathbf{d}^{T} \mathbf B^T \mathbf{\Phi_d})^{-1} \mathbf{\Psi}_\mathbf{d}^{T} \mathbf{L)V}
\end{aligned}
\end{equation}

\noindent and
\begin{equation}
    \mathbf{K_c} = \boldsymbol{T}_{\mathbf{d}}^T  \mathbf{K} \boldsymbol{T_d}  \in \mathbb{R}^{m\times m}.
\end{equation}

\noindent The very same operator $\boldsymbol{T}$ maps interface displacements to the subdomain via the relation
\begin{equation} \label{downscaling disp}
 u^h(x) = \mathbf N(x) \boldsymbol T \mathbf{u_{\Gamma}}.
\end{equation}

In this work, our objective is to parameterize the operators $\boldsymbol{T}$ and $ \mathbf{K_c}$ as functions of geometric parameters to lift the computational cost during the optimization loop.

\subsection{EIFEM-Based Preconditioning Strategy} \label{sec:EIFEM preconditioner}
In this section, we recall the EIFEM methodology, originally introduced in \cite{RUBIO2025} as a preconditioning strategy for linear systems. In the present work, however, within the optimization loop, a single EIFEM application (without pre- or post-smoothing) replaces the exact solution of the state PDE. Solving the PDE only approximately at each iteration yields, in turn, an inexact gradient of the objective function, a concept that is formalized in Section \ref{Section:optimization}.


As we discussed in \cite{RUBIO2025}, one could use EIFEM to obtain a fast approximation of a FE solution or one can use it as a coarse space within a preconditioner strategy. The benefit of the preconditioning strategy is that the full order FE solution is sought, while the computational time reduced.

In the context of solving large-scale systems of equations $\mathbf{Ax} = \mathbf{b}$ arising from PDE discretizations, preconditioners accelerate the convergence of the solver by reducing the condition number of the system, $\kappa(A)$. Let $N$ denote the total number of fine-scale DoFs in the global system, such that the system matrix is $\mathbf{A} \in \mathbb{R}^{N \times N}$ and the vectors are $\mathbf{x}, \mathbf{b} \in \mathbb{R}^{N}$.

Following the principles of domain decomposition \cite{TosselliWidlund2006,Doolean2015,Farhat2001FETIDP,mandel2008multispace} and multigrid methods \cite{Trottenberg2000,Olshanskii2014,Briggs2000}, an effective preconditioner must target both high- and low-frequency error components. To achieve this, we propose a multiplicative preconditioning scheme where EIFEM targets low-frequency errors while pre- and post-smoothing steps eliminate high-frequency errors.

The application of the preconditioner to the residual vector $\mathbf{r}_k \in \mathbb{R}^{N}$ at iteration $k$ is defined by the following three-step sequence, producing intermediate update vectors $\mathbf{z} \in \mathbb{R}^{N}$:
\begin{equation}
\setlength{\jot}{2pt}
    \begin{aligned}
    \mathbf{z}_{1/3} &= \mathbf{S}_{\nu_1} \mathbf{r}_k  \\
    \mathbf{z}_{2/3} &= \mathbf{z}_{1/3} + \boldsymbol{\phi}^T \left( \boldsymbol{\phi} \mathbf{A} \boldsymbol{\phi}^T \right)^{-1} \boldsymbol{\phi} (\mathbf{r}_k - \mathbf{A} z_{1/3})  \\
    \mathbf{z}_{3/3} &= \mathbf{z}_{2/3} + \mathbf{S}_{\nu_2} (\mathbf{r}_k - \mathbf{A} \mathbf{z}_{2/3}). 
\end{aligned}
\label{eq:pre_smooth}
\end{equation}

Here, the pre- and post-smoothing operators $\mathbf{S}_{\nu_1}, \mathbf{S}_{\nu_2} \in \mathbb{R}^{N \times N}$ represent $\nu_1$ and $\nu_2$ sweeps of a general smoothing matrix, respectively. Let $N_c$ denote the total number of DoFs in the coarse space ($N_c \ll N$). The operator $\boldsymbol{\phi} \in \mathbb{R}^{N_c \times N}$ represents the coarse-fine mapping, which projects the local EIFEM kinematics into the global fine-scale system. Note that the inverted term $\left( \boldsymbol{\phi} \mathbf{A} \boldsymbol{\phi}^T \right)^{-1} \in \mathbb{R}^{N_c \times N_c}$ represents the exact solve on the reduced EIFEM space.  

We define the global mapping operator $\phi$ explicitly as an assembly of $N_S$ subdomain contributions

\begin{equation} \label{eq:coarse_fine_mapping}
    \boldsymbol{\phi} = \sum_{i=1}^{N_S} \mathbf{\mathcal{A}}_i \boldsymbol{T}_i^T \mathbf{D}_i \mathbf{R}_i^T,
\end{equation}

\noindent where for each subdomain $\Omega_i$:
\begin{itemize}
    \item $\mathbf{R}_i \in \mathbb{R}^{N \times N_{loc}} $ is the restriction operator from the global DoFs to the local subdomain DoFs.
    \item $\mathbf{D}_i \in \mathbb{R}^{N_{loc} \times N_{loc}}$ is a diagonal scaling operator acting on the interface DoFs. It is defined by the counting operator $D_{jj} = 1/N_x$, where $N_x$ is the number of subdomains sharing node $j$.
    \item $\boldsymbol{T}_i \in \mathbb{R}^{N_{loc} \times m}$ is the interscale operator defined in Equation \eqref{downscaling disp}
    \item $\mathbf{\mathcal{A}}_i \in \mathbb{R}^{N_c \times m} $ is a coarse assembly operator.
\end{itemize}

As shown in \cite{RUBIO2025}, this formulation ensures that the proposed preconditioner is Symmetric Positive Definite (SPD) and suitable for the conjugate gradient method.

\section{Parameterized EIFEM} \label{Section: Parameterized EIFEM}

\subsection{Parametric PDEs}

In the context of parametric optimization, the physical system is not described by a single state, but by a collection of states governed by a set of parameters. We introduce a parameter domain $\mathcal{P} \subset \mathbb{R}^n$, representing the range of admissible configurations. A generic parameter vector $\boldsymbol{\mu} \in \mathcal{P}$ characterizes the variations in the system, which may include physical coefficients (e.g., conductivity, viscosity), boundary data, or geometric variations of the domain.

For a given parameter $\boldsymbol{\mu} \in \mathcal{P}$, we seek the field variable $u(\boldsymbol{\mu})$ that satisfies the partial differential equation in strong form:
\begin{equation}
    \mathcal{L}(u(\boldsymbol{\mu}); \boldsymbol{\mu}) = f(\boldsymbol{\mu}) \quad \text{in } \Omega(\boldsymbol{\mu}),
\end{equation}
subject to appropriate Dirichlet and/or Neumann boundary conditions on $\partial \Omega(\boldsymbol{\mu})$. Here, $\mathcal{L}(\cdot; \boldsymbol{\mu})$ represents a parameter-dependent differential operator of an elliptic problem, and $f(\boldsymbol{\mu})$ denotes the source term.

To allow for the application of numerical discretization techniques such as the FEM, we cast the problem in its variational form. Let $V(\boldsymbol{\mu})$ be an appropriate Hilbert space \big(e.g., the Sobolev space $H^1_0(\Omega(\boldsymbol{\mu}))$\big) defined over the domain. The weak formulation reads: find $u(\boldsymbol{\mu}) \in V(\boldsymbol{\mu})$ such that
\begin{equation}\label{eq:weak_form}
    a(u(\boldsymbol{\mu}), v; \boldsymbol{\mu}) = \ell(v; \boldsymbol{\mu}) \quad \forall v \in V(\boldsymbol{\mu}), 
\end{equation}
where $a(\cdot, \cdot; \boldsymbol{\mu}) : V(\boldsymbol{\mu}) \times V (\boldsymbol{\mu}) \to \mathbb{R}$ is a continuous bilinear form and $\ell(\cdot; \boldsymbol{\mu}) : V(\boldsymbol{\mu}) \to \mathbb{R}$ is a continuous linear form representing external loading.
We assume that for all $\boldsymbol{\mu} \in \mathcal{P}$, the bilinear form $a(\cdot, \cdot; \boldsymbol{\mu})$ is continuous and coercive ensuring the existence and uniqueness of the solution $u(\boldsymbol{\mu})$.

For efficient evaluation during iterative optimization, it is highly advantageous if the parameter dependence is affine. An operator is said to admit an affine decomposition if it can be expressed as:
\begin{equation} \label{affine decomp}
    a(u, v; \boldsymbol{\mu}) = \sum_{q=1}^{Q_a} \Theta_q^a(\boldsymbol{\mu}) \, a_q(u, v),
\end{equation}
where $\Theta_q^a(\boldsymbol{\mu})$ are a series of scalar functions depending \textit{only} on the parameter $\boldsymbol{\mu}$ and $a_q(u, v)$ are parameter-independent bilinear forms.
This structure allows the computationally expensive integrals (the matrices corresponding to $a_q$) to be pre-computed once (Offline), while the scalar parameter-dependent assembly (Online) becomes a trivial summation. The offline-online strategy is particularly appealing within an optimization loop since the matrices need to be constructed at every iteration, but this is done in the online phase. 

Although the affine decomposition in Equation (\ref{affine decomp}) is accurate in many physical problems, it is not always available, particularly in the presence of non-linearities or the case in which the geometry of the domain depends on the parameter. In such cases, the evaluation of the weak form integrals requires the assembly of the full high-fidelity operators at every iteration of the optimization process which negates the computational advantages of the parameterized model. To recover an efficient offline-online strategy, we use the Matrix version of the so-called Discrete Empirical Interpolation Discrete Empirical Interpolation Method (DEIM) \cite{Chaturantabut2010,quarteroni2016reduced,NEGRI2015}.

\subsection{The Discrete Empirical Interpolation Method (DEIM)}\label{DEIM section}
 Let $\mathbf{f}(\boldsymbol{\mu}) \in \mathbb{R}^{N_{loc}}$ denote the vector resulting from the discretization of a non-linear or non-affine term, for us the operators  $\boldsymbol{T}$ and $ \mathbf{K_c}$ defined in Equations \eqref{coarse problem} and \eqref{downscaling disp}. We seek to approximate this vector in a low-dimensional subspace spanned by a basis $\boldsymbol{W} = [\boldsymbol{w}_1, \dots, \boldsymbol{w}_k] \in \mathbb{R}^{N_{loc} \times k}$, where $k \ll N_{loc}$.
 
 The construction of $\boldsymbol{W}$ is performed in an offline training phase using the Proper Orthogonal Decomposition (POD). For this purpose, the parameter domain $\mathcal{P}$ is sampled at $S$ locations $\{\boldsymbol{\mu}_1, \dots, \boldsymbol{\mu}_S\}$ to generate a snapshot matrix $\mathbf{S}_f = [\mathbf{f}(\boldsymbol{\mu}_1), \dots, \mathbf{f}(\boldsymbol{\mu}_S)]$. By computing the Singular Value Decomposition (SVD), $\mathbf{S}_f = \mathbf{U} \mathbf{\Sigma} \mathbf{V}^T$, the basis $\boldsymbol{W}$ is chosen as the first $k$ columns of $\mathbf{U}$, corresponding to the dominant singular values. The approximation is then given by:

\begin{equation}
\mathbf{f}(\boldsymbol{\mu}) \approx \boldsymbol{W} \mathbf{c}(\boldsymbol{\mu}),
\end{equation}

\noindent where $\mathbf{c}(\boldsymbol{\mu}) \in \mathbb{R}^k$ is a vector of parameter-dependent coefficients. To determine $\mathbf{c}(\boldsymbol{\mu})$, DEIM selects a set of $k$ interpolation indices (the so-called "magic points"), denoted by $\mathcal{MP} = \{p_1, \dots, p_k\}$. We define a selection matrix $\mathbf{P} = [\mathbf{e}_{p_1}, \dots, \mathbf{e}_{p_k}] \in \mathbb{R}^{N_{loc} \times k}$, where $\mathbf{e}_{p_i}$ is the $p_i$-th column of the identity matrix. In simpler terms, $\mathbf{e}_{p_i}$ is a boolean vector that has only a non-zero entry at index $p_i$.

By enforcing the interpolation to be exact at these $k$ points, we obtain:
\begin{equation} \label{coefficientsDEIM}
\mathbf{P}^T \mathbf{f}(\boldsymbol{\mu}) = (\mathbf{P}^T \boldsymbol{W}) \mathbf{c}(\boldsymbol{\mu}).
\end{equation}

Provided the basis $\boldsymbol{W}$ and the indices $\mathcal{MP}$ are constructed such that the $k \times k$ matrix $(\mathbf{P}^T \boldsymbol{W})$ is non-singular, the approximation is given by:
\begin{equation}
\mathbf{f}(\boldsymbol{\mu}) \approx \boldsymbol{W} (\mathbf{P}^T \boldsymbol{W})^{-1} \mathbf{P}^T \mathbf{f}(\boldsymbol{\mu}).
\end{equation}

\noindent The computational efficiency of this approach stems from the fact that $\mathbf{P}^T \mathbf{f}(\boldsymbol{\mu})$ only requires evaluating the non-linear function at the $k$ sampled indices. Consequently, the online complexity becomes independent of the full mesh size $N_{loc}$, depending only on the reduced dimension $k$.

The selection of the interpolation indices $\mathcal{MP}$ is performed via a greedy strategy, detailed in Algorithm \ref{alg:DEIM}. The process initializes by selecting the index where the first basis vector has its maximum magnitude. In subsequent steps $\ell = 2, \dots, k$, the algorithm removes the contribution of the basis vectors already selected from the current mode $\boldsymbol{w}_\ell$. It then identifies the index $p_\ell$ where the remaining residual error $\mathbf{r}$ is maximal.  This ensures that the interpolation points are placed in regions of highest non-linear activity, minimizing the global approximation error in the max-norm.

\begin{algorithm}[H] 
\caption{Magic points selection algorithm \cite{Chaturantabut2010}}
\label{alg:DEIM}
\begin{algorithmic}[1]
\REQUIRE Basis vectors $\{\boldsymbol{w}_{\ell}\}_{\ell=1}^m \subset \mathbb{R}^{N_{loc}}$ obtained via POD on non-linear snapshots.
\ENSURE Indices $\mathcal{MP} = [p_1, \dots, p_k]^T$ and selection matrix $\mathbf{P}$.

\STATE \textbf{Initialization:}
\STATE Select the first index corresponding to the maximum magnitude of the first basis vector:
\begin{equation*}
p_1 = \arg\max_{i \in \{1, \dots, N_{loc}\}} |\boldsymbol{w}_1(i)|
\end{equation*}
\STATE Set $\mathbf{U} = [\boldsymbol{w}_1]$, $\mathbf{P} = [\mathbf{e}_{p_1}]$, and $\mathcal{MP} = [p_1]$.

\FOR{$\ell = 2$ to $k$}
    \STATE Solve the linear system for coefficients $\mathbf{c}$:
    \begin{equation*}
    (\mathbf{P}^T \mathbf{U}) \mathbf{c} = \mathbf{P}^T \boldsymbol{w}_{\ell}
    \end{equation*}
    \STATE Compute the residual (the error between the basis vector and its current approximation):
    \begin{equation*}
    \mathbf{r} = \boldsymbol{w}_{\ell} - \mathbf{U}\mathbf{c}
    \end{equation*}
    \STATE Select the next index where the residual error is largest:
    \begin{equation*}
    p_{\ell} = \arg\max_{i \in \{1, \dots, N_{loc}\}} |\mathbf{r}(i)|
    \end{equation*}
    \STATE Update the matrices:
    \STATE $\mathbf{U} \leftarrow [\mathbf{U}, \boldsymbol{w}_{\ell}]$
    \STATE $\mathbf{P} \leftarrow [\mathbf{P}, \mathbf{e}_{p_{\ell}}]$
    \STATE $\mathcal{MP} \leftarrow [\mathcal{MP}, p_{\ell}]$
\ENDFOR
\STATE \textbf{Output:} The set of indices $\mathcal{MP} = \{p_1, \dots, p_l\}$.
\end{algorithmic}
\end{algorithm}

\subsection{Parameterized EIFEM operators}
\begin{figure}[hbt!] 
\centering
\includegraphics[width=0.7\columnwidth]{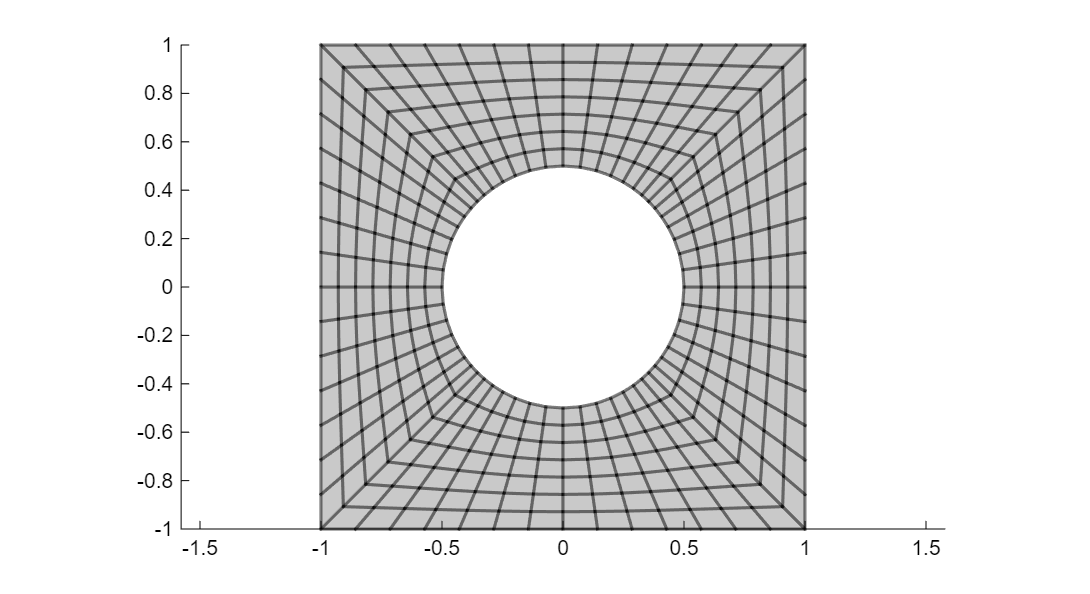}
\caption{FE mesh of the cell showing the circular inclusion. The parameter is the radius of the inclusion.}
\label{meshHole}
\end{figure}

As already mentioned, in this work we aim at accelerating the optimization process. To do so, we will make use of EIFEM to approximate both the energy functional and the state equation (see Section \ref{Parametric Optimization}). Therefore, we first need to parameterize the operators $\boldsymbol{T}(\boldsymbol{\mu})$ and $ \mathbf{K_c}(\boldsymbol{\mu})$ by means of the DEIM. Note that $ \mathbf{K_c}(\boldsymbol{\mu})$ could be constructed from $\boldsymbol{T}(\boldsymbol{\mu})$ via the relation $ \mathbf{K_c}(\boldsymbol{\mu}) = \boldsymbol{T}(\boldsymbol{\mu})^T \mathbf{K}(\boldsymbol{\mu}) \boldsymbol{T} (\boldsymbol{\mu})$ but we choose to approximate the coarse stiffness matrix $\mathbf{K_c}$ directly since it lies in a lower dimensional space than $\mathbf{K}$ and this fact drastically reduces the number of basis vectors resulting from DEIM.
To illustrate and validate the procedure, we will focus on a one-dimensional parameter space. We consider a square domain with a circular hole, where the parameter is the radius of the inclusion, as depicted in Figure \ref{meshHole}.

Following Section \ref{DEIM section}, we first sample the parametric space $\mathcal{P}$ at evenly spaced $S$ locations $\{\boldsymbol{\mu}_1, \dots, \boldsymbol{\mu}_S\}$ across the parameter range to generate the snapshot matrices $\mathbf{S}_T = [\boldsymbol{T}(\boldsymbol{\mu}_1), \dots, \boldsymbol{T}(\boldsymbol{\mu}_S)]$ and  $\mathbf{S}_\mathbf{K_c} = [\mathbf{K_c}(\boldsymbol{\mu}_1), \dots, \mathbf{K_c}(\boldsymbol{\mu}_S)]$, and then perform the SVD. Note that, both $\boldsymbol{T}  \in \mathbb{R}^{N_{loc} \times m}$ and $\mathbf{K_c} \in \mathbb{R}^{m \times m} $ are matrices, hence to store them we concatenate their columns in a single vector column, thus $\mathbf{S}_T \in \mathbb{R}^{(N_{loc} \cdot m ) \times S} $ and  $\mathbf{S}_\mathbf{K_c}  \in \mathbb{R}^{( m \cdot m ) \times S}$.

\begin{figure}[htbp]
    \centering
    \begin{subfigure}{0.65\textwidth}
        \includegraphics[width=\textwidth]{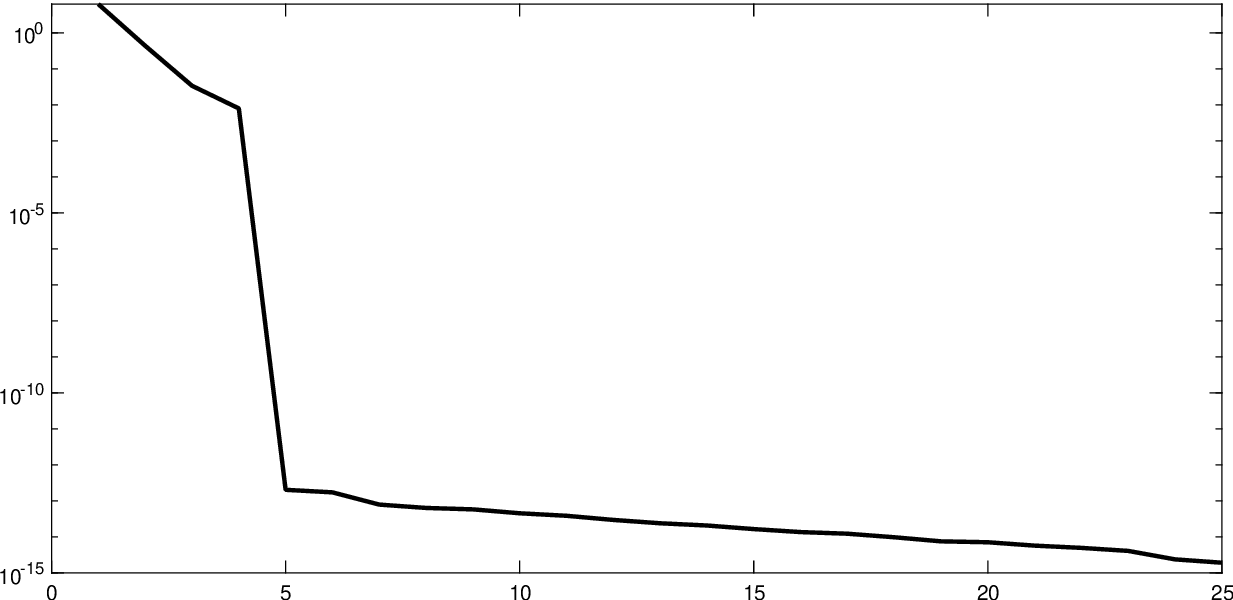}
        \caption{Decay of the singular values of $\mathbf{K_c}$. The presence of the spectral gap indicates that 4 modes approximate well the energy of the system.}
        \label{kSingular}
    \end{subfigure}
    \begin{subfigure}{0.65\textwidth}
        \includegraphics[width=\textwidth]{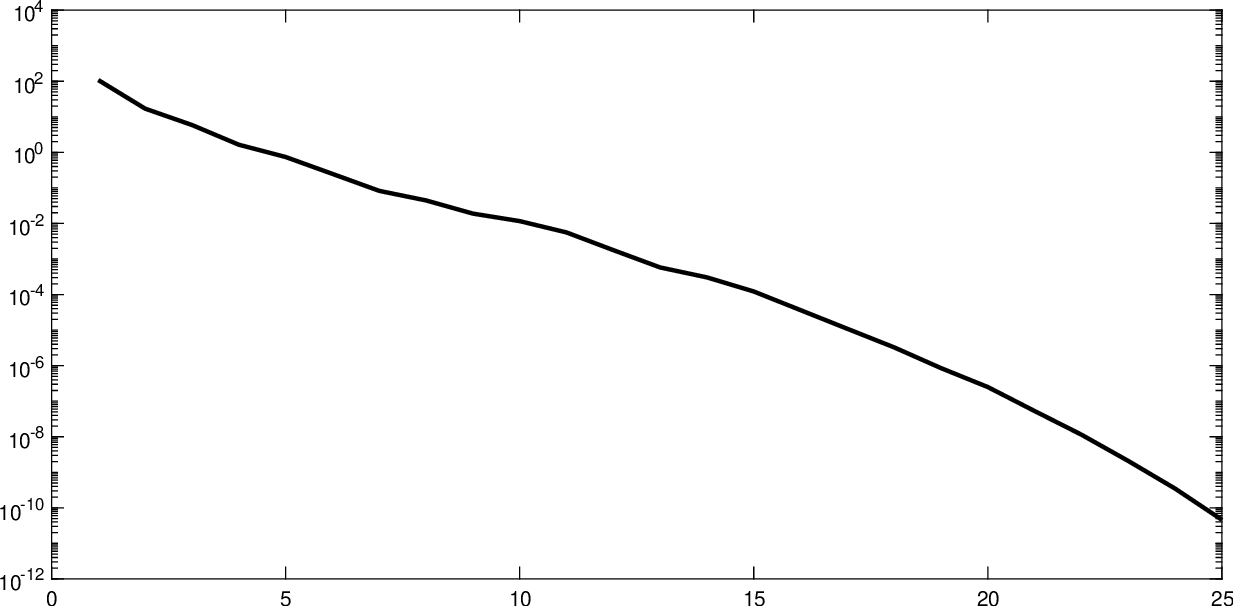}
        \caption{Decay of the singular values of $\boldsymbol{T}$. The absence of spectral gap suggests that the energy is more spread through the spectrum, however we observe an exponential decay of the singular values.}
        \label{TSingular}
    \end{subfigure}

    \caption{Decay of the singular values as a function of the mode index : (a) $\boldsymbol{K_c}$ and (b) $\boldsymbol{T}$ .}
    \label{fig:singularVal}
\end{figure}

Figure \ref{kSingular} presents the decay of the singular values of $\mathbf{K_c}$ (in logarithmic scale) as a function of their index. The spectrum shows two distinctive patterns. The first few singular values dominate the energy content of the system, with a relatively gradual decrease up to the fourth mode. This is followed by a sharp drop of several orders of magnitude, indicating a spectral gap. Such a rapid initial decay is characteristic of systems with low-dimensional structure suggesting that the dominant dynamics can be captured accurately by a small number of modes.

Figure \ref{TSingular} shows the singular value decay of $\boldsymbol{T}$ on a logarithmic scale as a function of the mode index. In contrast to the previous case, the spectrum presents a smooth decay without a pronounced spectral gap. The singular values decrease approximately exponentially (nearly linear trend in the semi-log representation), spanning several orders of magnitude across the full set of modes.

The absence of a clear separation between dominant and non-dominant modes indicates that the energy of the system is distributed more broadly across the basis. Although the leading modes still carry the largest contribution, higher-order modes may retain non-negligible energy, and their decay is progressive rather than abrupt. This behavior suggests a more complex or less strongly compressible structure compared to the previous spectrum.

\begin{figure}[hbt!] 
\centering
\includegraphics[width=0.65\columnwidth]{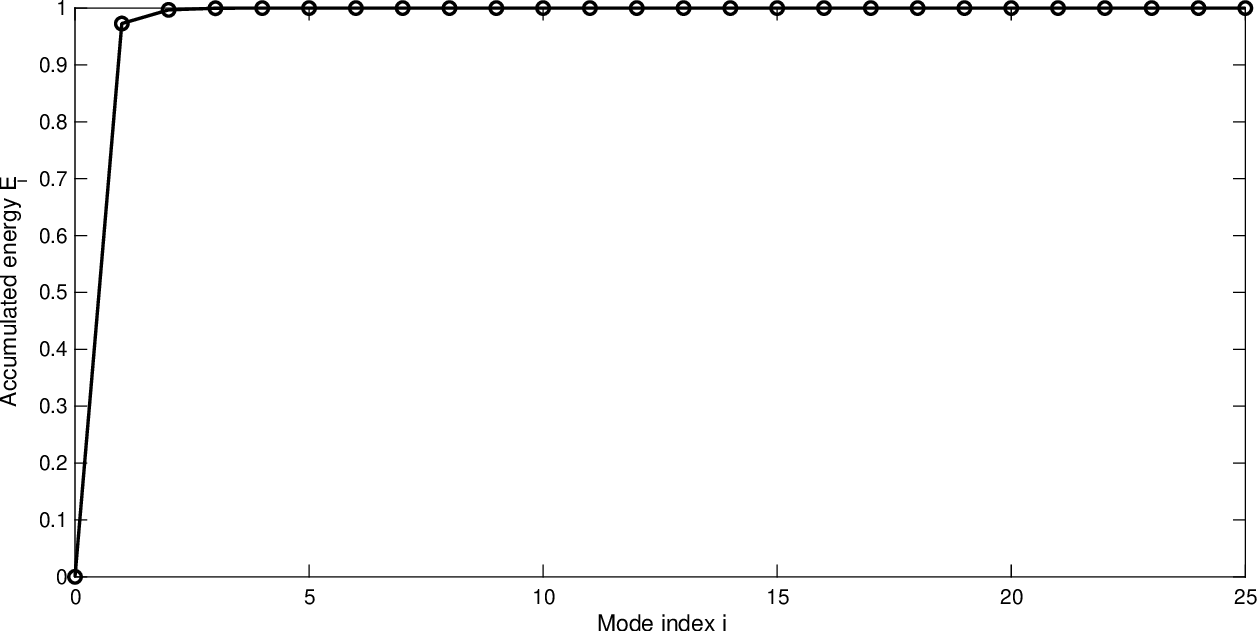}
\caption{Accumulated energy across the spectra of $\boldsymbol{T}$. We observe that the leading modes are clearly dominant which indicates that few of the can accurately reproduce the original operator. Setting the error tolerance $\epsilon_k = 1-E_K$ to $1 \cdot 10^{-5}$, the algorithm retrieves 10 modes. }
\label{Tenergy}
\end{figure}

To quantify the truncation efficiency, Figure \ref{Tenergy} shows the cumulative energy per mode index, defined as

\begin{equation}
E_k = \frac{\sum_{i=1}^k \sigma_i^2}{\sum_{j=1}^{S} \sigma_j^2},  
\end{equation}

\noindent where $\sigma_i$ denotes the i-th singular value. Despite the smoother decay of the individual singular values, the accumulated energy rises rapidly, capturing most of it within the first modes. This confirms that a reduced basis is still viable, though it may require a larger dimension $k$ compared to $\mathbf{K_c}$ to achieve the same accuracy. 

The cumulative energy thus allows us to define a criterion on the number of basis vectors to retrieve. Of course, there is a trade-off between the approximation error and the number of basis vectors. In this case, since the number of basis vectors is small, we can push the error tolerance $\epsilon_k = 1-E_K$ to $1 \cdot 10^{-5}$, which results in 10 modes for $\boldsymbol{T}$.

\begin{figure}[hbt!] 
\centering
\includegraphics[width=1.05\columnwidth]{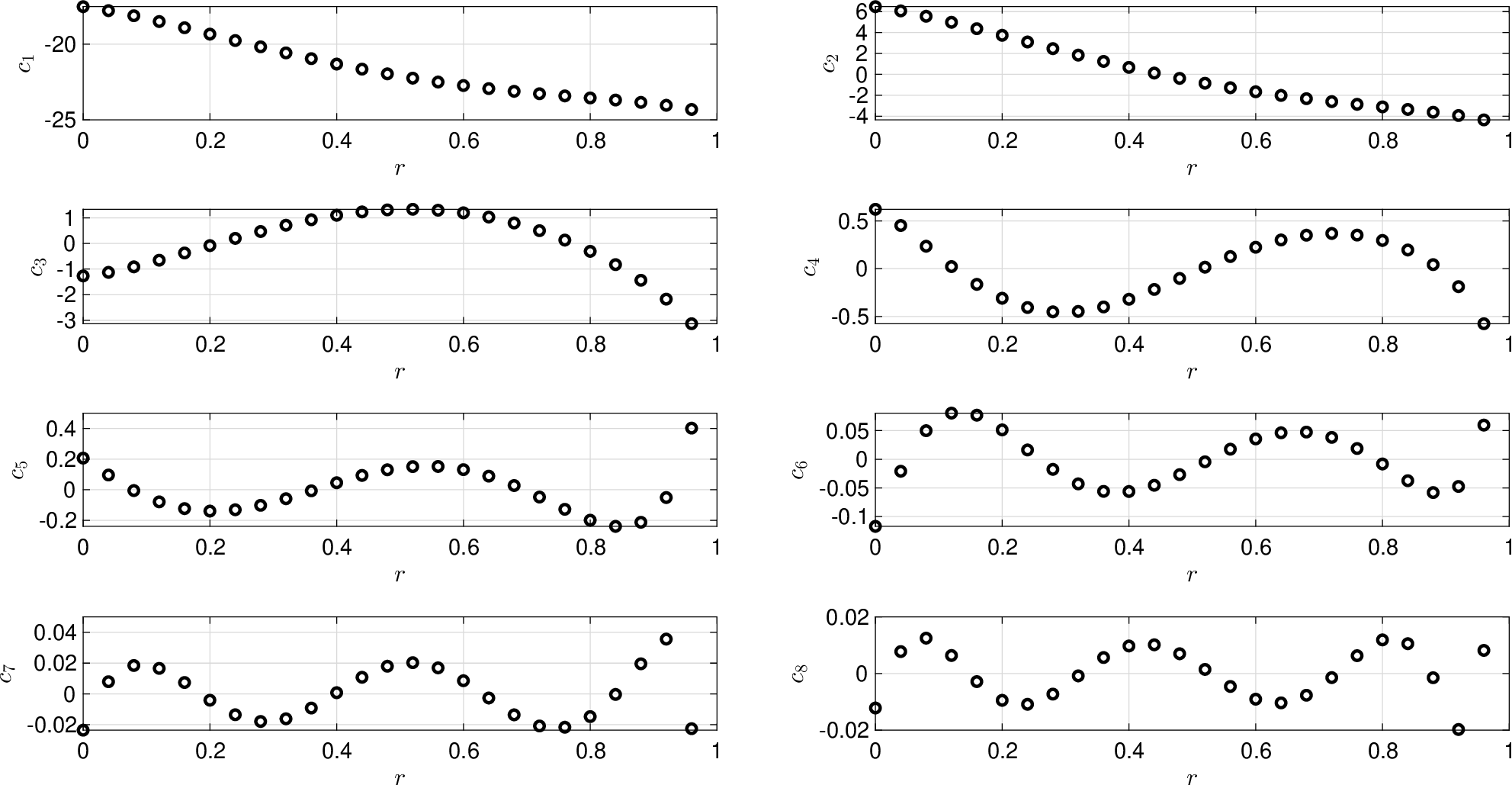}
\caption{Evolution of the first 6 coefficients as a function of the radius for $\boldsymbol{T}$.}
\label{coefficients_T}
\end{figure}

After defining the criterion to select $\boldsymbol{W} = [\boldsymbol{w}_1, \dots, \boldsymbol{w}_k]$, we can run Algorithm \ref{alg:DEIM} to obtain the indices $\mathcal{MP} = \{p_1, \dots, p_k\}$ and, once the indices are obtained, we can apply Equation (\ref{coefficientsDEIM}) to compute the dependency of the coefficients $\mathbf{c}(\boldsymbol{\mu})$ with respect to the parameters, as shown in Figure \ref{coefficients_T}.

At this point, the vector of coefficients $\mathbf{c}(\boldsymbol{\mu})$ is available only at the discrete parameter samples $\{\boldsymbol{\mu}_1, \dots, \boldsymbol{\mu}_S\}$, but we seek a continuous approximation over the parameter domain $\boldsymbol{\mu} \in \mathcal{P}$. Given the smooth dependence of the coefficients on $\boldsymbol{\mu}$, we decided to approximate the data with a Lagrange polynomial of degree $\polDeg=8$. For each component of the coefficient vector $\mathbf{c}_i(\boldsymbol{\mu})$, the approximation reads:
\begin{equation}
\tilde{\mathbf{c}}_i(\boldsymbol{\mu}) = \sum_{j=1}^{\polDeg+1} \alpha_{ij} \phi_j(\boldsymbol{\mu}), \quad \text{for } i = 1, \dots, k,
\label{eq:poly_approx}
\end{equation}
where $\phi_j(\boldsymbol{\mu})$ are the Lagrange basis polynomials and $\alpha_{ij}$ are the coefficients to be determined via least squares minimization. For every coefficient we solve

\begin{equation}
\left\{
\begin{aligned}
\min_{y,\alpha} \quad & \frac{1}{2}\|y-c(\boldsymbol{\mu})\|_2^2 \\[-1mm]
\text{s.t.} \quad & y = Z \alpha
\end{aligned}  
\right.
\end{equation}

\noindent where $ Z \in  \mathbb{R}^{s \times (\polDeg+1)} $ is the matrix that contains the evaluation of the shape functions $\phi_j(\boldsymbol{\mu})$ at the sampling point $\{\boldsymbol{\mu}_1, \dots, \boldsymbol{\mu}_S\}$. 

An important detail is that, when high-degree Lagrange polynomials are used, evenly spaced nodes lead to oscillations near the interval boundaries, known as Runge’s phenomenon \cite{Canuto2006}, \cite{Hesthaven_2007}. To avoid this issue, non-uniform nodes such as Chebyshev, Gauss–Lobatto, or Gauss–Legendre points should be employed. In this work we used Chebyshev points. 

The continuous mapping $\boldsymbol{W} \tilde{\mathbf{c}}_i(\boldsymbol{\mu})$ allows the operators $\boldsymbol{T}(\boldsymbol{\mu})$ and $ \mathbf{K_c}(\boldsymbol{\mu})$ to be quickly reconstructed for any new radius $\boldsymbol{\mu}$ using straightforward matrix vector products and scalar evaluations. This step completely decouples the online evaluation phase from the full-order dimension, giving the desired computational speed-up.

Before moving into the validation phase, let us comment on how to extend the interpolation to higher dimensions. A natural extension to higher-dimensional parameter spaces, i.e., $\boldsymbol{\mu} \in \mathbb{R}^n$ with $n > 1$, can be obtained by tensor-product of one-dimensional Lagrange polynomials defined along each parameter direction. However, this strategy suffers from the so-called curse of dimensionality, i.e., the number of basis vectors scales exponentially as $(\polDeg+1)^n$. To mitigate the exponential growth of degrees of freedom associated with tensor grids, one may employ, for example, sparse grids \cite{PfluegerP2010}.

\subsection{Validation results: EIFEM vs parametric EIFEM}
The validation of the method is performed in two steps.  First, we compute the approximation error across the parameter space 
$\mathcal{P}$  for both $\boldsymbol{T}(\boldsymbol{\mu})$ and $ \mathbf{K_c}(\boldsymbol{\mu})$ using the Frobenius norm
\begin{equation} \label{FrobeniusNorm}
\frac{\|A(\boldsymbol{\mu}) - \widetilde{A}(\boldsymbol{\mu})\|_F}
     {\|A(\boldsymbol{\mu})\|_F}
= \frac{\sqrt{\sum_{i=1}^{m} \sum_{j=1}^{n} 
\left| a_{ij}(\boldsymbol{\mu}) - \widetilde{a}_{ij}(\boldsymbol{\boldsymbol{\mu}}) \right|^2}}{\sqrt{\sum_{i=1}^{m} \sum_{j=1}^{n} 
\left| a_{ij}(\boldsymbol{\mu}) \right|^2}},
\end{equation}

\noindent where $A(\boldsymbol{\mu})$ is the exact and $\tilde{A}(\boldsymbol{\mu})$ the approximated operator. Second we assess the influence of the parameterized operators on the convergence of the preconditioner described in Section \ref{sec:EIFEM preconditioner}. 

\begin{figure}[htbp]
    \centering
    \begin{subfigure}{0.65\textwidth}
        \includegraphics[width=\textwidth]{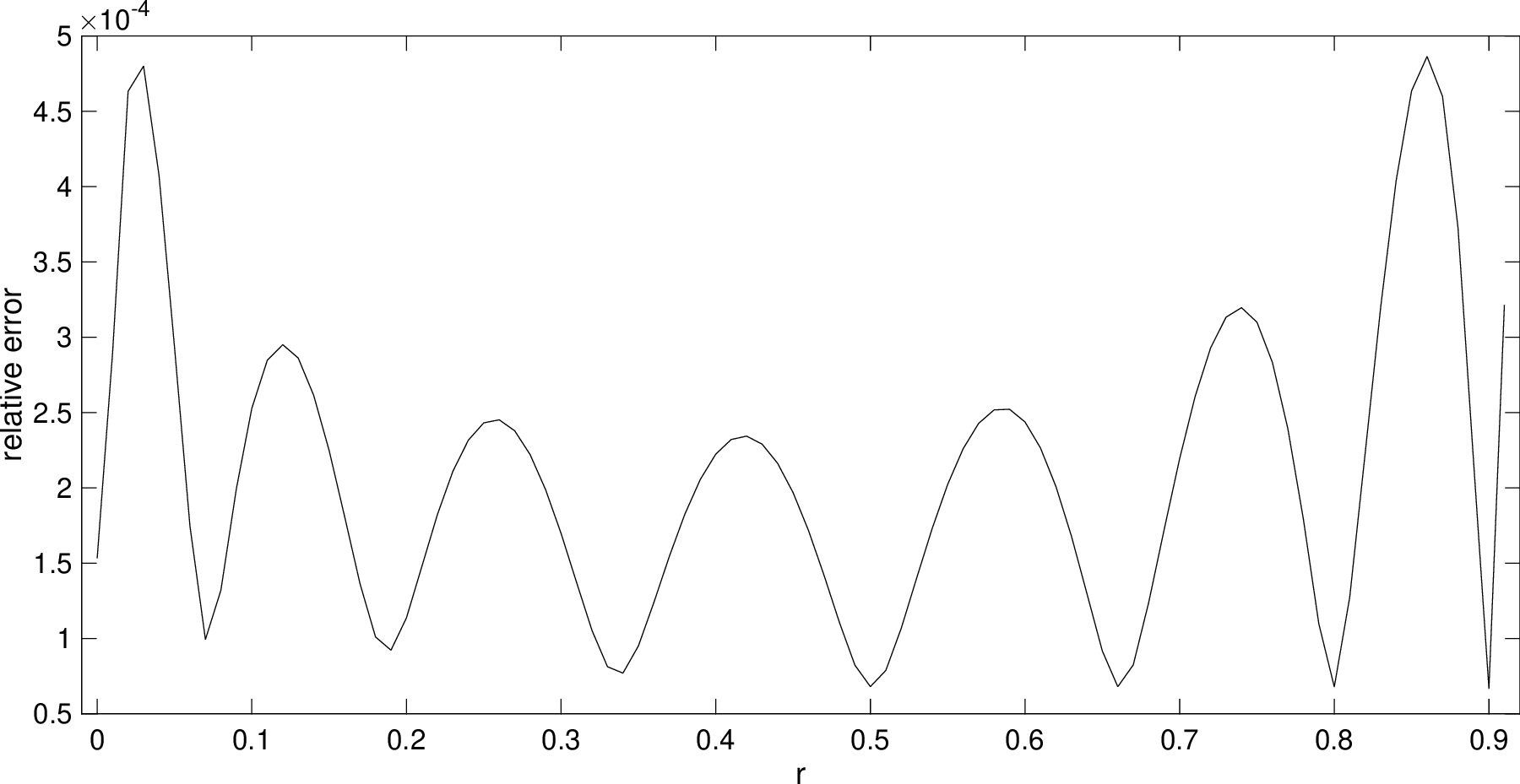}
        \caption{$\boldsymbol{T}(\boldsymbol{\mu})$.}
        \label{error_T}
    \end{subfigure}
    \quad
    \begin{subfigure}{0.65\textwidth}
        \includegraphics[width=\textwidth]{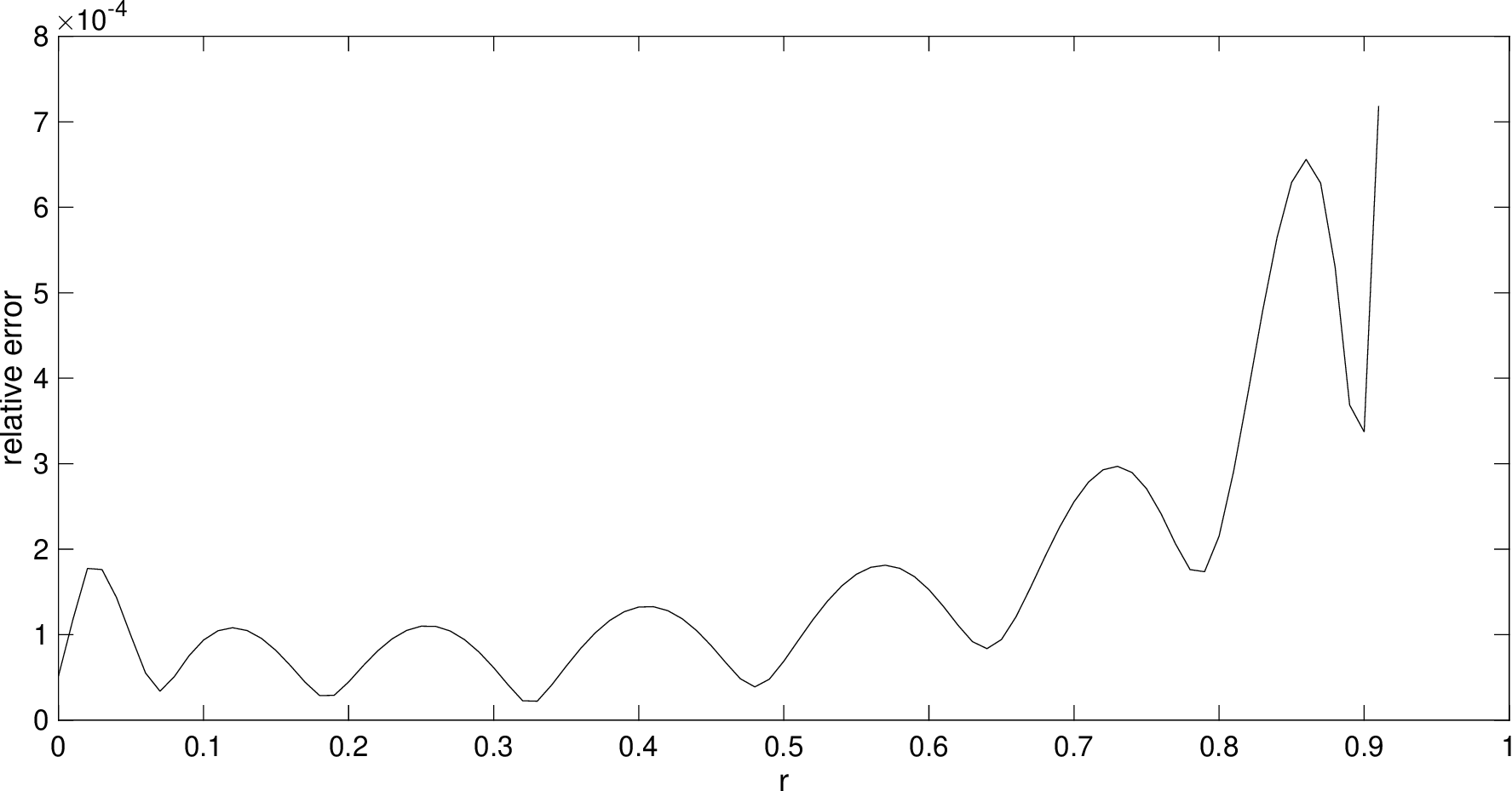}
        \caption{$\boldsymbol{K_c}(\boldsymbol{\mu})$.}
        \label{error_Kc}
    \end{subfigure}

    \caption{Approximation error as a function of the parameter $r$: (a) $\boldsymbol{T}(\boldsymbol{\mu})$ and (b) $\boldsymbol{K_c}(\boldsymbol{\mu})$.}
    \label{fig:error_approx}
\end{figure}

Figure \ref{error_T} shows that the relative error associated with $\boldsymbol{T}(\boldsymbol{\mu})$ remains uniformly small throughout the parameter domain,  with values on the order of $10^{-4}$. The error has a smooth oscillatory behaviour but no significant growth, 
indicating a stable and accurate approximation across the parameter space.

Figure \ref{error_Kc} presents the relative error for $ \mathbf{K_c}(\boldsymbol{\mu})$. The error is also of the order $10^{-4}$ for most values of the parameter. Although a moderate increase is observed for larger values of $r$, the approximation remains accurate throughout the domain.

Overall, the results demonstrate that the parametric EIFEM provides a highly accurate approximation of both operators. The relative error consistently remains below $10^{-3}$ and is typically of order $10^{-4}$ over the whole parameter space, confirming the robustness and reliability of the proposed approach.

\begin{figure}[hbt!] 
\centering
\includegraphics[width=0.85\columnwidth]{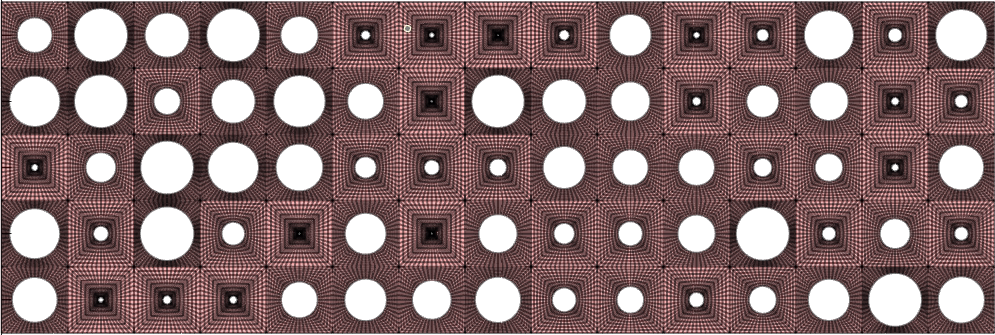}
\caption{FE mesh of randomly distributed radius that is used for the validation of the interpolated operators.}
\label{validationMesh}
\end{figure}

The second stage of the validation is the impact of the parameterized operators in the convergence of the preconditioner compared to exact ones. Figure \ref{validationMesh} shows the FE mesh of a randomly distributed parameter r. The parameter values are assigned randomly to the subdomains to generate a heterogeneous test case, allowing the performance and convergence of the preconditioner based on parameterized operators to be assessed under non-uniform subdomain meshes.

Table \ref{tab:exact vs interpolated iters} compares the performance of the interpolated and exact operators for a random distribution of the parameter $r$ over the range $10^{-6}-0.96$. The results show that the number of iterations remains identical for both approaches, indicating that the interpolated operators preserve the effectiveness of the preconditioner. Furthermore, reducing the number of snapshots from 100 to 25 does not affect convergence, showing that a relatively small sampling set is sufficient to capture the parametric dependence. In terms of computational cost, the interpolated approach reduces the solution time from approximately 4 minutes to about 10 seconds, yielding more than an order of magnitude speedup. These results highlight that the proposed interpolation strategy achieves substantial computational savings while maintaining the efficiency of the exact operators.

\begin{table}[!h]
\centering
\begin{tabular}{|p{0.14\textwidth}|p{0.14\textwidth}|p{0.14\textwidth}|p{0.10\textwidth}|p{0.10\textwidth}|}
\hline 
 & \# Snapshots & \# DEIM Basis & Iterations & Time \\
\hline 

\multirow{3}{0.14\textwidth}{Interpolated}
& 100 & 10 & 248 &  4.5 sec \\
\cline{2-5}

& 50 & 10  & 248 & 4.5 sec \\
\cline{2-5}

& 25 & 10 & 248 & 4.5 sec \\
\hline

Exact & $\emptyset$ & $\emptyset$ & 248 & 4 min \\
\hline

\end{tabular}
\caption{Comparison of the performance of the interpolated and exact operators for a random distribution of the parameter r ranging from $1e-6$ to $0.96$. The number of domains for all cases is $45 \times15$. The iteration counts for both approaches are identical which validates the procedure followed. As for the number of snapshots, the results indicate that 25 sampling points suffice to have a good performance. By exact operators we mean that every subdomain has been trained; that is why the running time is much higher than for the interpolated ones. The interpolated times only include solve time.}
\label{tab:exact vs interpolated iters}
\end{table}

\section{EIFEM-based structural optimization} \label{Section:optimization}
Before restricting the design space to a set of geometric parameters, let us first briefly introduce the problem of TopOpt, since this will help us present how parametric optimization can be helpful to speed up the optimization process.  

\subsection{Topology Optimization} \label{sec:topology_optimization}
Let the material domain $\Omega $ be a subset of a larger reference domain $ D \subset\mathbb{R}^d$. We denote the  boundary of $D$ by $\partial D$. The characteristic function $\chi:D \rightarrow \{0,1\}$ defined as
\begin{equation}
    \chi(x) = 
    \begin{cases} 
        1 & \text{if } x \in \Omega \\ 
        0 & \text{if } x \notin \Omega 
    \end{cases},
    \label{eq:chi_def}
\end{equation}

\noindent indicates the presence of solid material at a given point in the reference domain. $x \in D$ denotes the position vector. 

In TopOpt one seeks to minimize a given functional $J(\chi,u)$ subject to some PDE constraints and, possibly, a set of additional constraints $g(\chi)$:
\begin{equation}
    \begin{cases}
        \min\limits_{\chi,u} & J(\chi, u) \\
        \text{s.t.} & a(\chi, u, v) = \ell(v) \quad \forall v \in H^1_0(D) \\
        & g(\chi) \le 0
    \end{cases},
    \label{eq:to_general}
\end{equation}

\noindent where $a(\chi, u, v)$ is the bilinear form that depends on the characteristic function $\chi$, $\ell(v)$ the linear form and $u \in H^1(D)$ is the solution of the PDE \cite{BendsoeSigmund2003,Sigmund2013}. However, note that (\ref{eq:to_general}) is not differentiable with respect to $\chi$, therefore, in order to have differentiability, the problem is approximated by means of density \cite{BendsoeSigmund2003} or level-set functions \cite{Allaire2004}. Let us just mention the former since it will serve as smooth transition between the original TopOpt problem (\ref{eq:to_general}) and the parametric optimization we employ.

In the density approach, the problem is relaxed and, rather that using the characteristic function $\chi$, a continuous version of it namely density $\rho: D \rightarrow [0, 1]$ is introduced. Then

\begin{equation}
    \begin{cases}
        \min\limits_{\rho,u} & J(\rho, u) \\
        \text{s.t.} & a(\rho, u, v) = \ell(v) \quad \forall v \in H^1_0(D) \\
        & g(\rho) \le 0
    \end{cases},
    \label{eq:to_density}
\end{equation}

\noindent which, particularized for structural optimization with linear elasticity and volume and box constraints for $\rho$ reads:

\begin{equation}
    \begin{cases}
        \min\limits_{\rho,u} & \ell(u) \\
        \text{s.t.} & a(\rho, u, v) = \ell(v) \quad \forall  v \in H^1_0(D) \\
        &  \int_{D} \rho(x) \, d\Omega - V_f \int_{D} 1 \, d\Omega \leq 0 \\
        & 0 \leq \rho(x) \leq 1 \quad \forall x \in D \\
    \end{cases},
    \label{eq:compliace_to_density}
\end{equation}

\noindent where $V_f$ is the target volume fraction and

\begin{equation}
    a(\rho, u, v) = \int_D \nabla^S u : \mathbb{C}(\rho) : \nabla^S v \,  dV , 
    \label{eq:bilinear_elasticity}
\end{equation}

\begin{equation}
    \ell(v) = \int_D vb  \, dV + \int_{\partial D} v t \, dS .
    \label{eq:inear_elasticity}
\end{equation}

\noindent where $b$ denotes the body forces and $t$ the tractions.

To solve the optimization problem efficiently, the constrained problem is transformed into an unconstrained (or strictly bound-constrained) form. The core idea is to eliminate the PDE constraint and embed it directly into the objective function.
 
Assuming well-posedness of the bilinear form $a(\rho, u, v)$ for any admissible distribution $\rho$, there exists a unique $u(\rho) \in H^1(D) $ solving  $a(\rho,u(\rho),v) = l(v) \; \forall v$. Therefore 

\begin{equation}
    \begin{cases}
        \min\limits_{\rho} & \ell(u(\rho)) \\
        \text{s.t.} 
        &  \int_{D} \rho(x) \, d\Omega - V_f \int_{D} 1 \, d\Omega \leq 0 \\
        & 0 \leq \rho(x) \leq 1 \quad \forall x \in D \\
    \end{cases},
    \label{eq:compliace_to_density_reduced}
\end{equation}

\subsubsection{Sensitivity analysis}
The reduced formulation is solved via gradient based algorithms, as the Method of Moving Asymptotes (MMA) \cite{Svanberg1987MMA} or the Null space algorithm \cite{Feppon2020} which require the evaluation of the functional as a function of the design variable as well as its derivative. In practice, the continuous problem described in Equation (\ref{eq:compliace_to_density_reduced}) is discretized using the FEM leading to

\begin{equation}
    \begin{cases}
        \min\limits_{\rho} & \mathbf{F}^T \mathbf{u(\rho)} \\
        \text{s.t.} 
        &  \int_{D} \rho(x) \, d\Omega - V_f \int_{D} 1 \, d\Omega \leq 0 \\
        & 0 \leq \rho(x) \leq 1 \quad \forall x \in D \\
    \end{cases}.
    \label{eq:compliace_to_density_reduced_discrete}
\end{equation}

\noindent where $\mathbf{u}(\rho) \in \mathbb{R}^N$ solves the linear system $\mathbf{K}(\rho) \mathbf{u} = \mathbf{F}$.

The gradient of the cost functional with respect to $\rho$ can then be computed via chain rule as

\begin{equation}
   \frac{dJ}{d\rho} = \frac{\partial J}{\partial \mathbf{u}}\frac{\partial \mathbf{u}}{\partial \rho} = \mathbf{F}^T \mathbf{K}^{-1} \left( -\frac{\partial  \mathbf{K}}{\partial \rho} \mathbf{u} \right) =  - \mathbf{u}^T \frac{\partial  \mathbf{K}}{\partial \rho}  \mathbf{u} \ .
    \label{eq:gradient}
\end{equation}

\noindent Since the problem is self-adjoint, evaluating the gradient requires solving $\mathbf{u}(\rho)$ (not the adjoint), which we recall is typically the computational bottleneck of the algorithm.

\subsection{Parametric Optimization} \label{Parametric Optimization}
In order to speed up the optimization process, we leverage three main aspects. First, the PDE will be approximated with EIFEM. Second, the derivative of the coarse stiffness matrix is computed analytically, since the operators $\mathbf{K_c}$ and $\boldsymbol{T}$ are explicitly parameterized with DEIM and high-order Lagrange polynomials (see Section \ref{Section: Parameterized EIFEM}). Finally, instead of discretizing $\rho$ using FE, the design variables $\boldsymbol{\mu}$ are introduced as geometric parameters associated with each subdomain, reducing the dimensionality of the optimization problem. In this way, a coarsening strategy is applied to both the design and the state variables.

More precisely, let the domain $D$ be decomposed into $N_S$ non-overlapping subdomains such that
\begin{equation}
    D = \bigcup_{i=1}^{N_S} D_i,
\qquad
D_i \cap D_j = \emptyset \ \text{for } i \neq j.
\end{equation}

\noindent The design variable is defined at the subdomain level as 
\begin{equation}
    \boldsymbol{\mu} = (\boldsymbol{\mu}_1, \dots, \boldsymbol{\mu}_{N_s}) \in \mathbb{R}^{N_S},
\end{equation}
hence, the number of design variables (which can be more than $1$ per subdomain) depends on the number of subdomains and is independent of the finite element discretization. We can then approximate problem (\ref{eq:compliace_to_density_reduced_discrete}) as

\begin{equation}
\begin{cases}
\displaystyle \min_{\boldsymbol{\mu},\mathbf{u_c}} 
& \mathbf{\tilde{F}_c}^T \mathbf{u_c} \\[2mm]
\text{s.t.} 
& \mathbf{\tilde K_c}(\boldsymbol{\mu}) \, \mathbf{u_c} = \mathbf{\tilde F_c}  \\[2mm]
& \ g(\boldsymbol{\mu}) - V_f |D| \le 0  \\[2mm]
& \boldsymbol{\mu}_{\min} \le \boldsymbol{\mu} \leq \boldsymbol{\mu}_{\max} .
\end{cases}
\end{equation}

\noindent where $g(\boldsymbol{\mu}) = \sum^{N_S}_{i=1} g_i(\boldsymbol{\mu_i})$ is an explicit function that returns the volume of the domain as a function of the parameters, while $\mathbf{\tilde F_c} $ denotes the coarse scale forces. In the most general case, $\mathbf{\tilde F_c}$ depends on the parameter via the relation $\mathbf{\tilde F_c}(\boldsymbol{\mu}) = \boldsymbol{\tilde T}(\boldsymbol{\mu})^T \mathbf{F}$. However, consistent with standard topology optimization formulations, we assume design-independent loading. Because the fine-scale external load $\mathbf{F}$ is assumed to remain constant regardless of the topology, we adopt the same strategy for the coarse-scale forces. Therefore, the variation of the compliance functional is only driven by the coarse displacement field $\mathbf{u_c}(\boldsymbol{\mu})$.

\subsubsection{Sensitivity analysis of the coarse problem}
Following the same reasoning as for the TopOpt case, we can define the reduced coarse functional

\begin{equation}
\begin{cases}
\displaystyle \min_{\boldsymbol{\mu},\mathbf{u_c}} 
& \mathbf{\tilde {F}_c}^T \mathbf{u_c}(\boldsymbol{\mu}) \\[2mm]
\text{s.t.} 
& \ g(\boldsymbol{\mu}) - V_f |D| \le 0  \\[2mm]
& \boldsymbol{\mu}_{\min} \le \boldsymbol{\mu} \leq \boldsymbol{\mu}_{\max} ,
\end{cases}
\label{coarse functional}
\end{equation}

\noindent where $\mathbf{u_c}(\boldsymbol{\mu}) \in \mathbb{R}^{N_c}$ ($N_c \ll N$) is the solution of the coarse problem $ \mathbf{\tilde K_c}(\boldsymbol{\mu}) \mathbf{u_c} = \mathbf{\tilde F_c}$. The gradient of the cost functional then becomes

\begin{equation}
   \frac{dJ}{d\boldsymbol{\mu}} = \frac{\partial J}{\partial \mathbf{u_c}}\frac{\partial \mathbf{u_c}}{\partial \boldsymbol{\mu}} = - \mathbf{u_c}^T \frac{\partial  \mathbf{\tilde K_c}}{\partial \boldsymbol{\mu}}  \mathbf{u_c} \ .
    \label{eq:gradient_coarse}
\end{equation}

Note that, since $\mathbf{\tilde K_c}$ is explicitly parameterized, its gradient only requires the evaluation of the derivatives of the coefficients with respect to the parameter

\begin{equation}
   \frac{d\mathbf{ \tilde K_c}}{d\boldsymbol{\mu}} =  \frac{d (\mathbf{W}\tilde{\mathbf{c}}(\boldsymbol{\mu}))}{d\boldsymbol{\mu}} = \mathbf{W} \frac{d  (\tilde{\mathbf{c}}(\boldsymbol{\mu}))}{d\boldsymbol{\mu}}
    \label{eq:gradient_coarse_param}
\end{equation}

\noindent instead of the more costly operation

\begin{equation}
   \frac{d\mathbf{K_c}}{d\boldsymbol{\mu}} =  \frac{\partial\boldsymbol{T}^T(\boldsymbol{\mu})}{ \partial \boldsymbol{\mu}} \mathbf{K}(\boldsymbol{\mu})\boldsymbol{T}(\boldsymbol{\mu}) + \boldsymbol{T}(\boldsymbol{\mu})^T  \frac{\partial\mathbf{K}(\boldsymbol{\mu})}{ \partial \boldsymbol{\mu}} \boldsymbol{T}(\boldsymbol{\mu}) +  \boldsymbol{T}(\boldsymbol{\mu})^T\mathbf{K}(\boldsymbol{\mu})\frac{\partial\boldsymbol{T}(\boldsymbol{\mu})}{ \partial \boldsymbol{\mu}} .
    \label{eq:gradient_coarse_param_full}
\end{equation}

\section{Numerical results} \label{Results}
\begin{figure}[hbt!] 
\centering
\includegraphics[width=1.05\columnwidth]{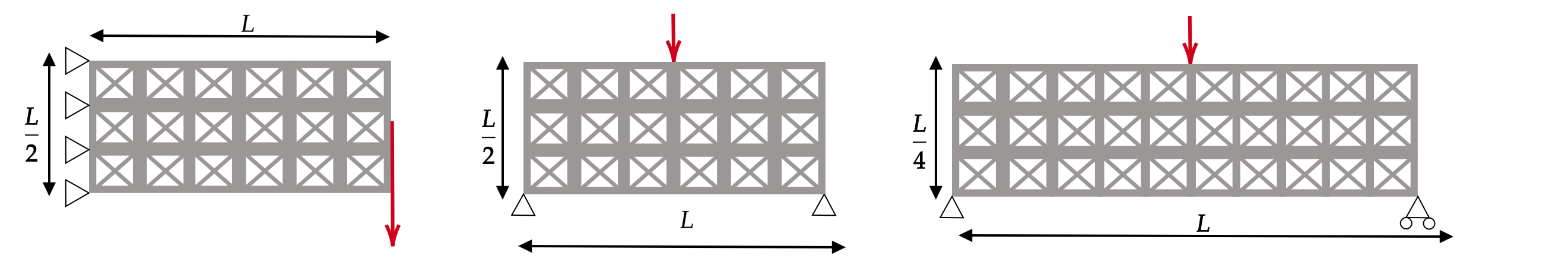}
\caption{Sketch of the benchmarks showing the boundary conditions. Cantilever beam (left), arch (middle) and MBB (right). }
\label{benchmarks}
\end{figure}

\begin{figure}[hbt!] 
\centering
\includegraphics[width=1\columnwidth]{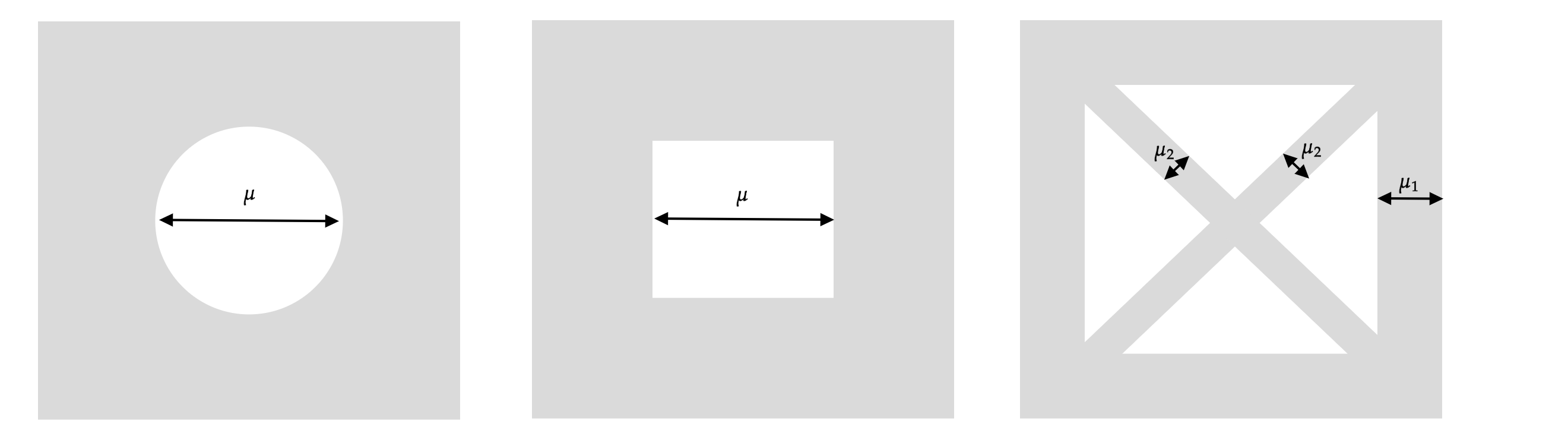}
\caption{Sketch of the geometry of unit cell of the circular inclusion (left), square inclusion (middle) and lattice cell (right) showing the geometrical parameters. }
\label{unitCells}
\end{figure}

In this section, the proposed parametric optimization strategy is assessed through a series of benchmark problems. The performance, accuracy, and computational efficiency of the parametric EIFEM approach are evaluated by comparing against TopOpt. The robustness of the algorithm is demonstrated across three classic structural optimization benchmarks: the Cantilever beam, the Arc, and the MBB beam (see Figure \ref{benchmarks}); and three different geometries: a circular inclusion, a square inclusion, and a lattice structure. The former two consist of a single geometrical parameter, while the geometry of the latter is controlled by two parameters, the width of the frames and the width of the bars as depicted in Figure \ref{unitCells}. We will use the the MMA \cite{Svanberg1987MMA} optimizer in all cases.

For the assessment of the designs obtained with parametric EIFEM, we will compare the final shape, the volume fraction distribution, the value of the compliance, and the computational time for all unit cells and TopOpt. In all cases, the EIF element considered is linear, consisting in 8 DoFs per subdomain in 2D and 24 DoFs in 3D

\subsection{Cantilever}
For the cantilever beam, we optimize a structure made of $28 \times 15 $ domains, which leads to a design space of 420 variables for single parameter geometries and 856 for the lattice cell. The parameter ranges from $10^{-6}$ to $0.96$ for the radius and square inclusion and between $0.06$ and $0.4$ for case of the lattice. The number of DoFs in the FE discretization is 330.000, however, we recall that in the parametric optimization we perform, the FE discretization is only needed for visualization purposes, all the operations are carried out at the coarse level. The volume target is $50\%$. For comparison with TopOpt, the number of degrees of freedom was reduced to 19,600, since higher resolutions result in substantially increased computational cost \footnote{The same discretization with 330,000 DoFs requires more than 8 hours to run on our machines. Therefore, we selected a mesh that is sufficiently fine while maintaining a reasonable computational cost.}. We stress that the focus of this study is on computational speed-up rather than on achieving the lowest possible compliance. This further highlights the appeal of EIFEM-based parametric operators: due to the coarsening strategy, large-scale optimization problems can be solved on a standard PC within seconds.

\begin{figure}[hbt!] 
\centering
\includegraphics[width=0.88\columnwidth]{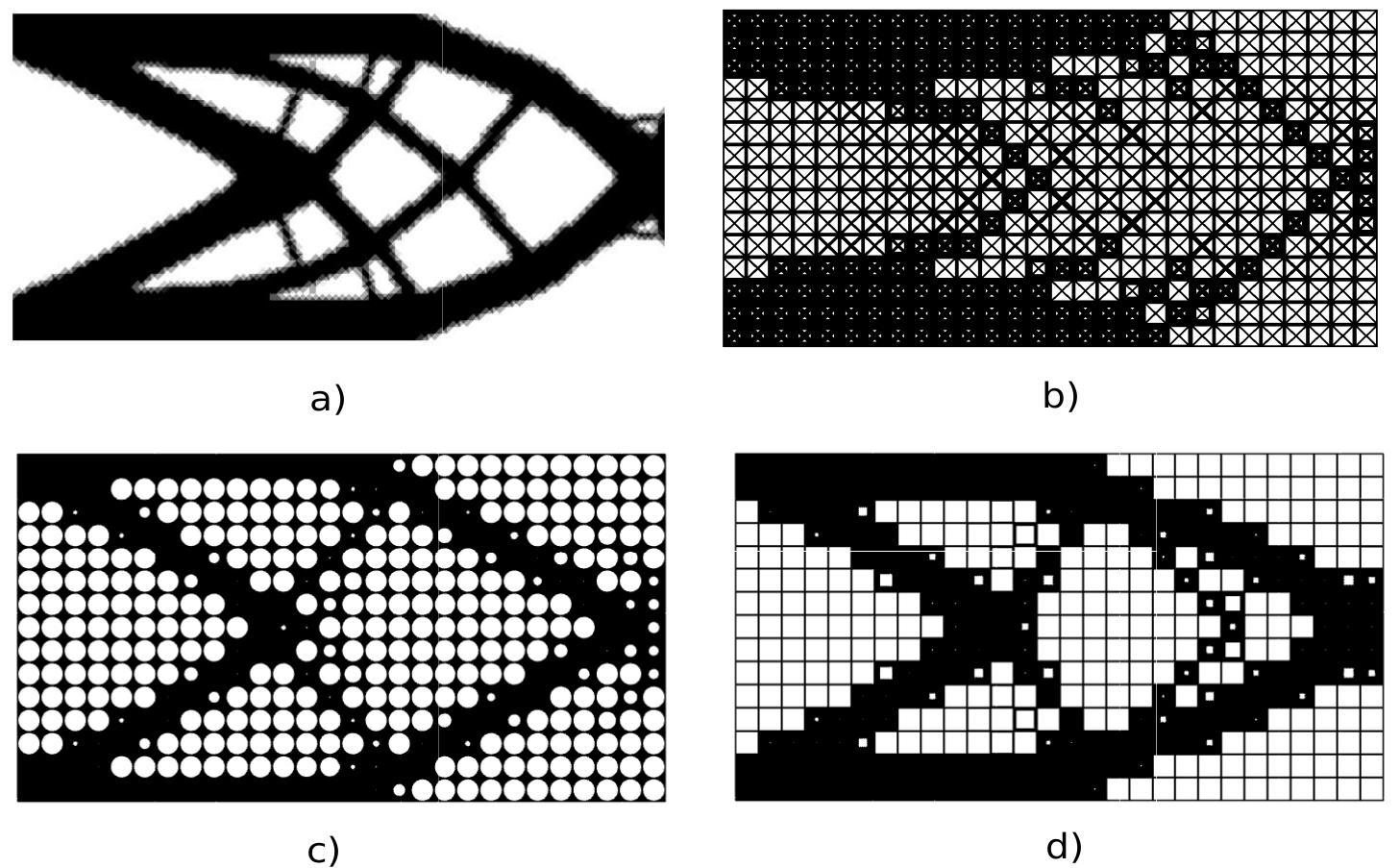}
\caption{Cantilever beam design obtained with TopOpt (a), the lattice unit cell (b) and the circular (c) and square (d) inclusions . A qualitatively similar structural pattern to classical TopOpt is observed. However, regions with intermediate (“gray”) densities appear when using the parametric optimization, especially for the lattice unit cell. In density-based TopOpt, such regions are typically penalized as they represent artificial mixtures of solid and void. In the present parametric framework, these regions correspond to continuous variations of geometric parameters and therefore emerge naturally when advantageous.}
\label{cantileverTopologies}
\end{figure}

In Figure \ref{cantileverTopologies}, we show that the cantilever beam design obtained from the lattice unit cell exhibits a qualitatively similar pattern to the one produced by classical TopOpt. However, regions with intermediate densities (“gray” areas) are observed.

In density-based TopOpt formulations, intermediate densities are typically penalized (e.g., via SIMP interpolation) because they represent mixtures of solid and void phases. Although such regions are often structurally competitive from a compliance minimization standpoint, they pose limitations in manufacturability.

\begin{figure}[hbt!] 
\centering
\includegraphics[width=0.88\columnwidth]{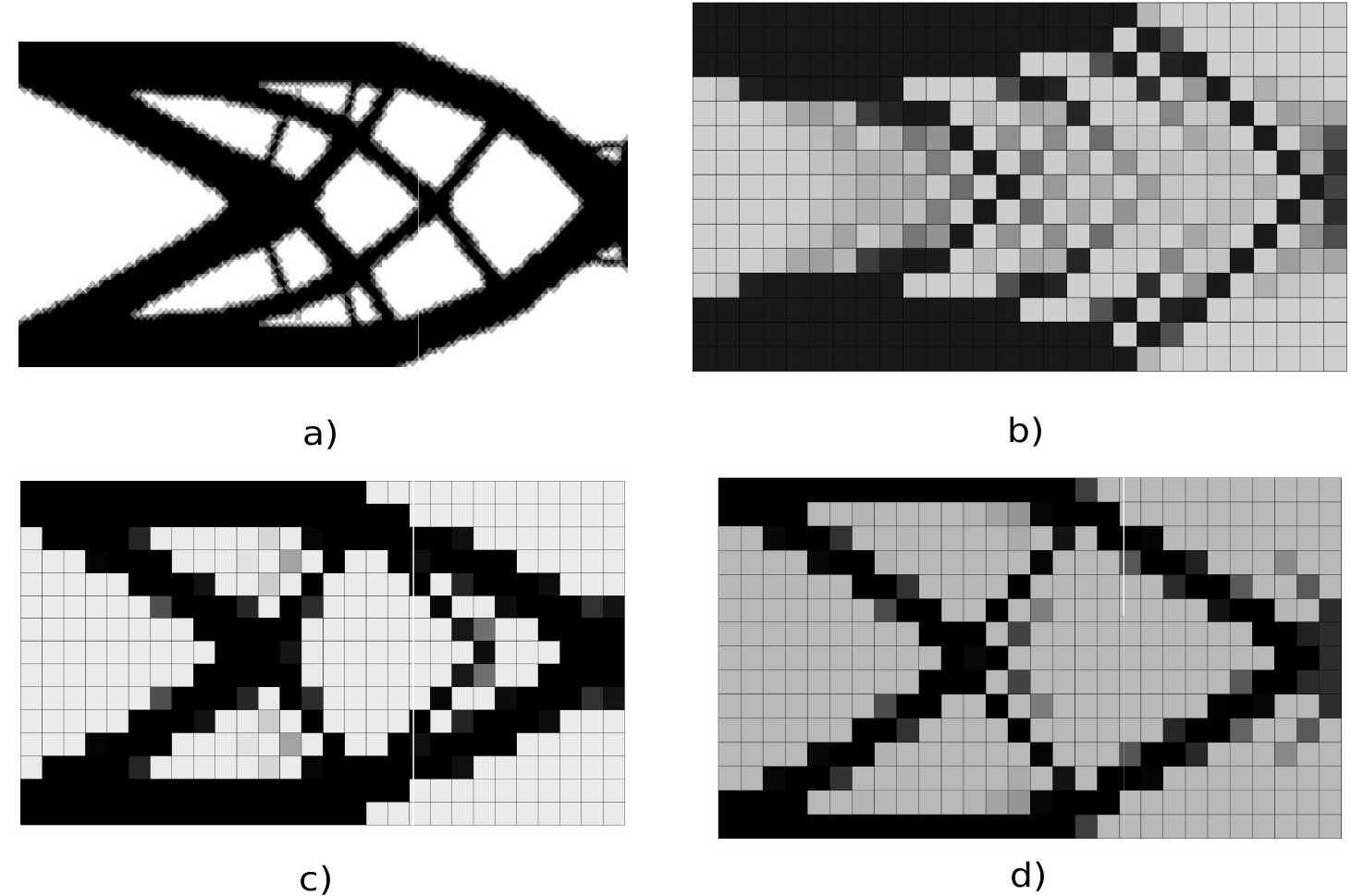}
\caption{Cantilever beam design obtained with TopOpt (a), the lattice unit cell (b) and the square (c) and circular (d) inclusions. A qualitatively similar structural pattern to classical TopOpt is observed. However, regions with intermediate (“gray”) densities appear when using the parametric optimization, especially for the lattice unit cell. In density-based TopOpt, such regions are typically penalized as they represent artificial mixtures of solid and void. In the present parametric framework, these regions correspond to continuous variations of geometric parameters and therefore emerge naturally when advantageous}.
\label{CantileverDensities}
\end{figure}
\begin{figure}[hbt!] 
\centering
\includegraphics[width=0.88\columnwidth]{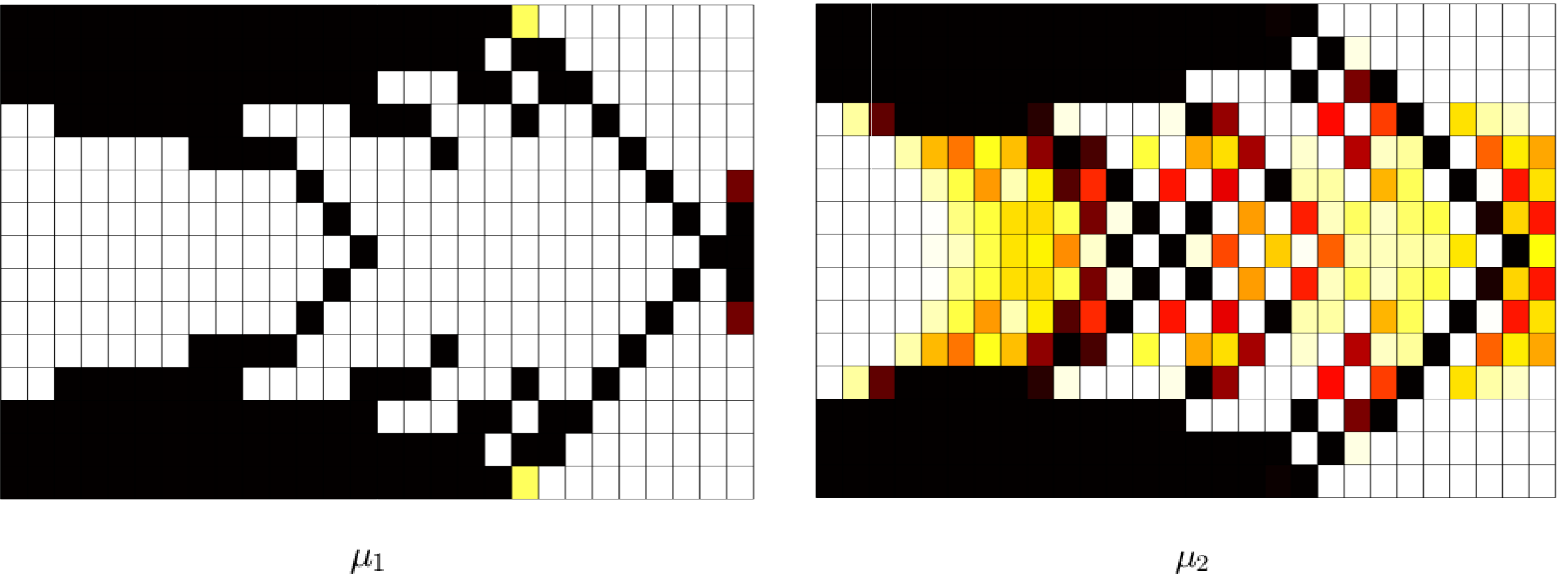}
\caption{Parameter distribution for the lattice cell. Black maximum value, white minimum value. In this configuration, the cross bars are aligned with the principal stresses, which promotes the inclusion of intermediate values of $\mu_2$ and is key for attaining stiff structures.}
\label{CantileverParameters}
\end{figure}

In contrast, with the parametric optimization framework adopted here, these intermediate regions arise naturally. Rather than representing fictitious material mixtures, they correspond to continuous variations of geometric parameters of the lattice microstructure. As a result, no artificial penalization is required, and geometrically “gray-like” regions can emerge as optimal solutions when they are advantageous (see Figures \ref{CantileverDensities}b and \ref{CantileverParameters}).

The circular inclusion produces a final design (Figure \ref{cantileverTopologies}c) that closely resembles the classical TopOpt baseline. Compared to the lattice unit cell case, the size of the “gray” regions is reduced, with a stronger tendency toward a black-and-white layout.

This behavior can be attributed to the geometry of the subdomain, which does not introduce preferential directions at the subdomain level. As a result, material is retained along the main load paths while removed from regions carrying lower loads. A similar mechanism is present in classical density-based TopOpt, where penalization promotes the formation of distinct solid and void regions. In the present framework, however, voids are introduced at the macroscale through geometric modification rather than through density interpolation. 

Similarly, the square inclusion leads to black-and-white structural features (Figure \ref{cantileverTopologies}d). This effect is amplified by the geometry of the inclusion, since it allows the removal of nearly all material from the cell, therefore facilitating sharper transitions between solid and void regions. Something worth mentioning is that stress concentrations near the corners may appear, but stress concentrations are not considered in this work as they do not affect the final geometry.

\begin{figure}[htbp]
    \centering
    \begin{subfigure}{0.58\textwidth}
        \includegraphics[width=\textwidth]{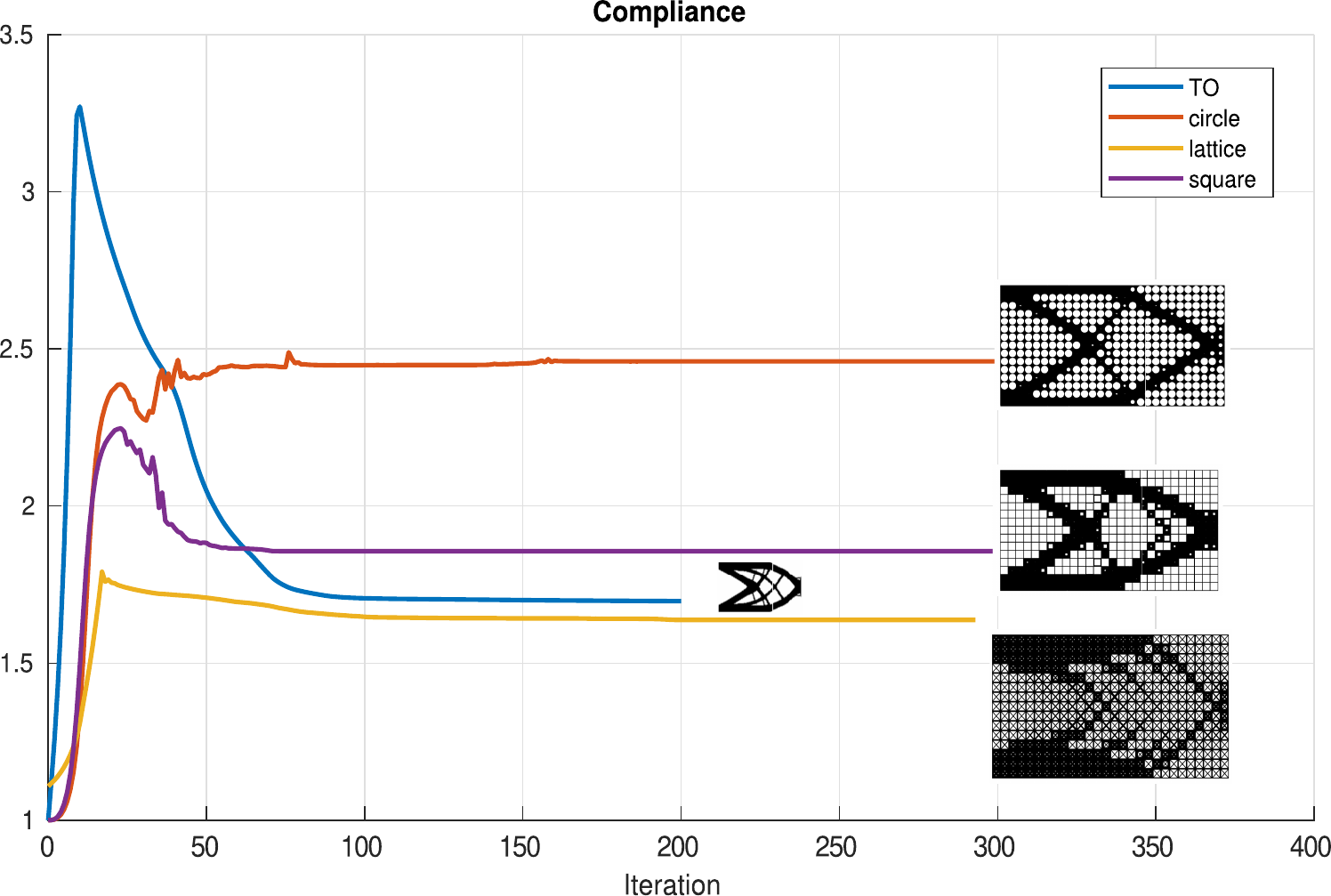}
        \caption{Convergence plots of the compliance objective function for the cantilever beam example. Results highlight the superior stiffness obtained with the lattice cell compared to classical TopOpt. Within the single parametric geometries, the square inclusion outperforms the circular inclusion.}
        \label{monitoringCantilever}
    \end{subfigure}
    \begin{subfigure}{0.72\textwidth}
        \includegraphics[width=\textwidth]{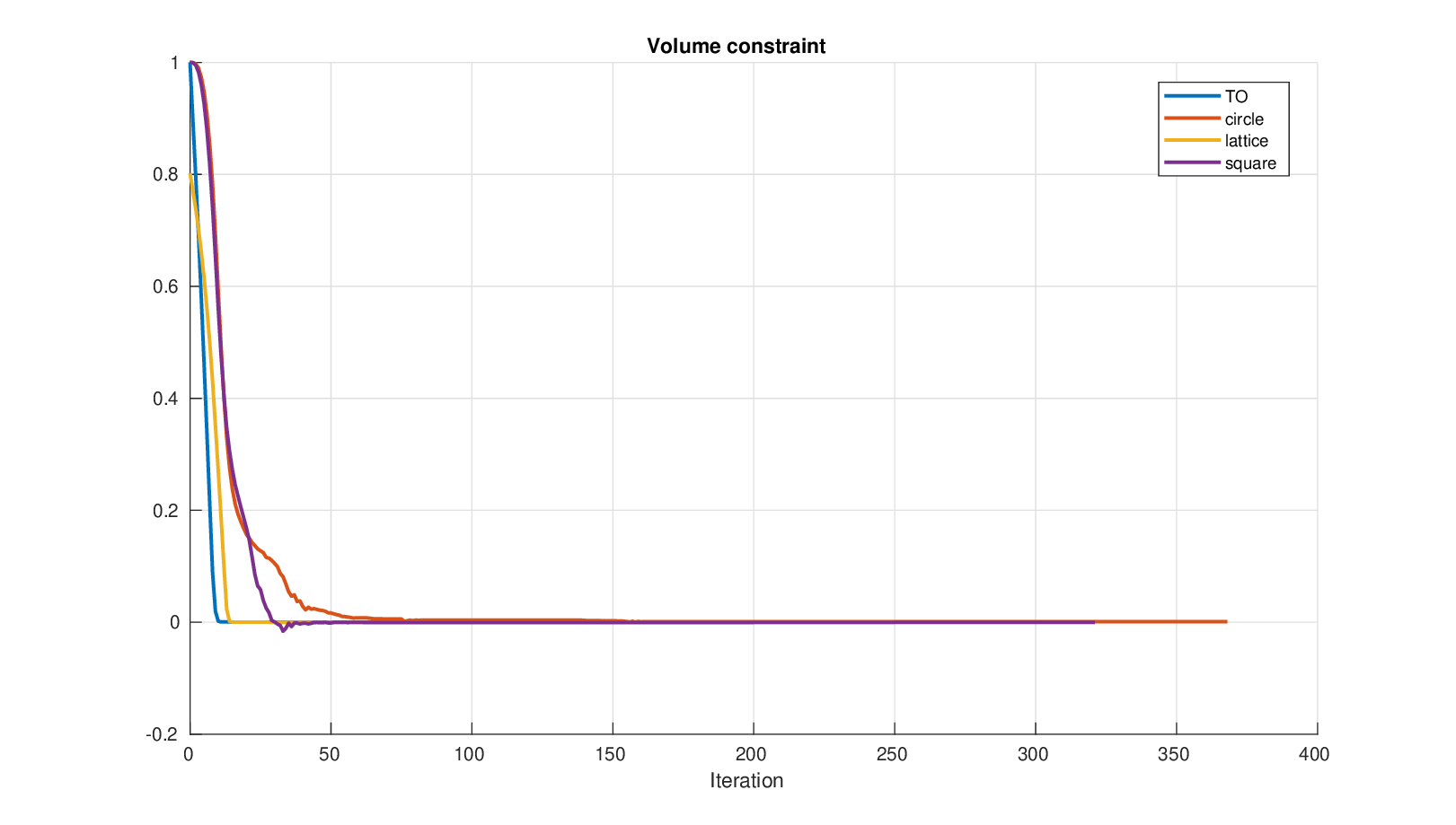}
        \caption{Relative volume constraint evolution as a function of the iteration.}
        \label{volumeCantilever}
    \end{subfigure}
    
    \caption{ Cantilever results. (a) Compliance evolution for TopOpt and the 3 different unit cells and (b) Volume evolution as a function of the iterations}
    \label{Cantilever results}
\end{figure}

To assess the performance of the designs, Figure \ref{monitoringCantilever} shows the evolution of the compliance functional during the optimization process. The results indicate that the circular inclusion leads to less efficient structures, whereas the lattice-cell parameterization yields the lowest compliance.

These two facts can be attributed to two main reasons. First, unlike classical topology optimization which penalizes intermediate densities, the lattice approach exploits these "gray" regions as a manufacturable unit cells. Because the lattice cell can transmit forces through both its frame and internal diagonal bars, its anisotropy promotes highly efficient, intermediate-density configurations. It is precisely this mechanism that allows the lattice design to achieve a globally stiffer structure than the classical TopOpt algorithm for this specific load case.

Second, when evaluating the other microstructures, performance is largely dictated by their void ratios. Parameterizations capable of removing a larger fraction of material from the unit cell provide the optimizer greater flexibility. This allows the mass to be redistributed more effectively to the regions where it  contributes to the structural stiffness. The total volume is the same for all geometries and TopOpt, as depicted in Figure \ref{volumeCantilever}.

\clearpage

\subsection{Arch}

In this second example, we use the same geometry and identical parameter range to those of the cantilever beam and the same discretization for the TopOpt baseline. The volume target is $50\%$.

\begin{figure}[hbt!] 
\centering
\includegraphics[width=0.88\columnwidth]{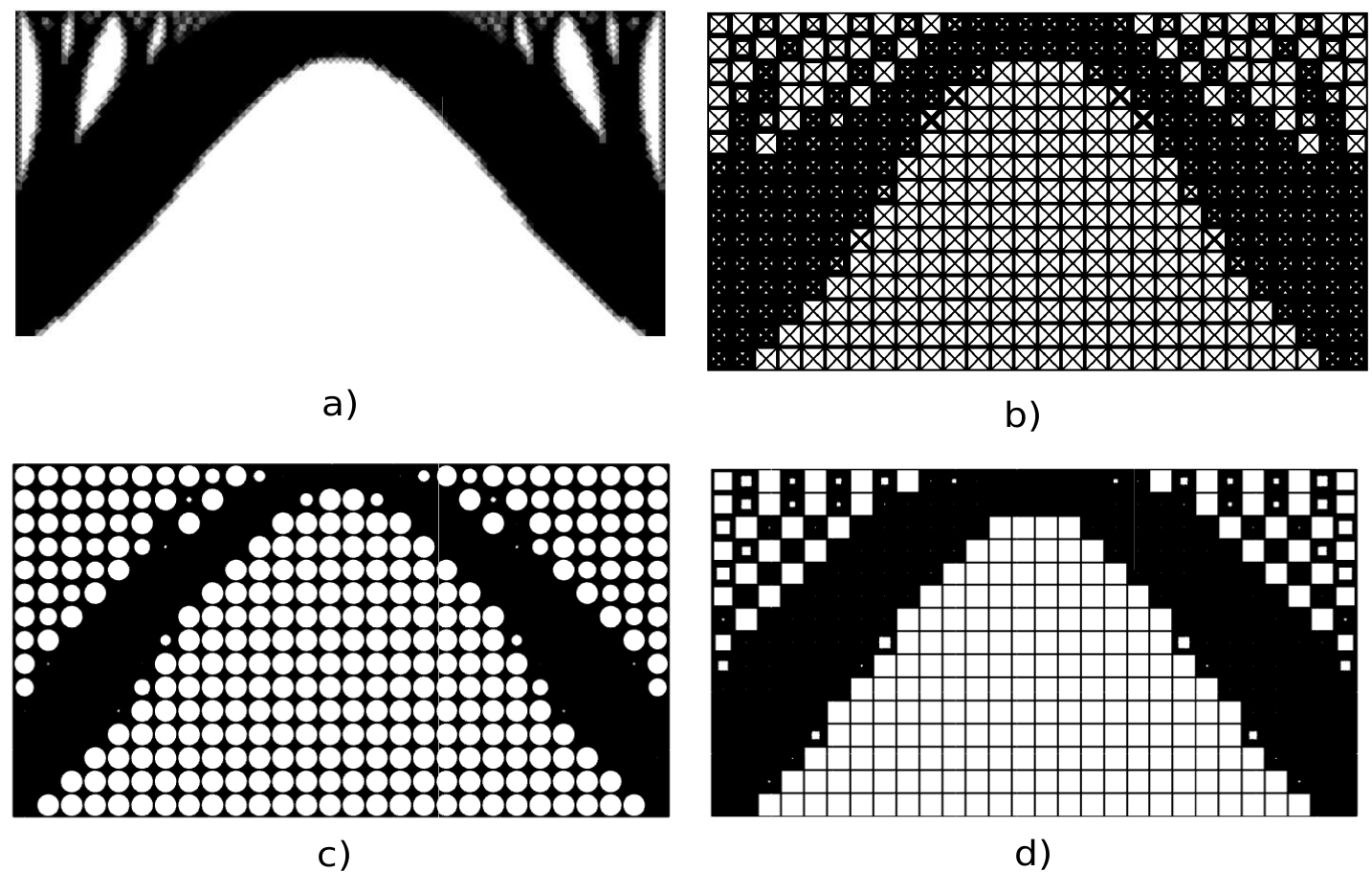}
\caption{Arch design obtained with a) TopOpt b) lattice unit cell c) circular and d) square inclusion. Parametric optimization obtains qualitatively similar structural pattern to classical TopOpt. In contrast to the cantilever beam, regions with intermediate densities are largely reduced, especially for the lattice unit cell. This is a consequence of the cross-bars being not aligned with the principal stresses and thus the only way to stiffen the structure is by adding material isotropically (thickening both the frame and the cross-bars).}
\label{archTopolgies}
\end{figure}

\begin{figure}[hbt!] 
\centering
\includegraphics[width=0.88\columnwidth]{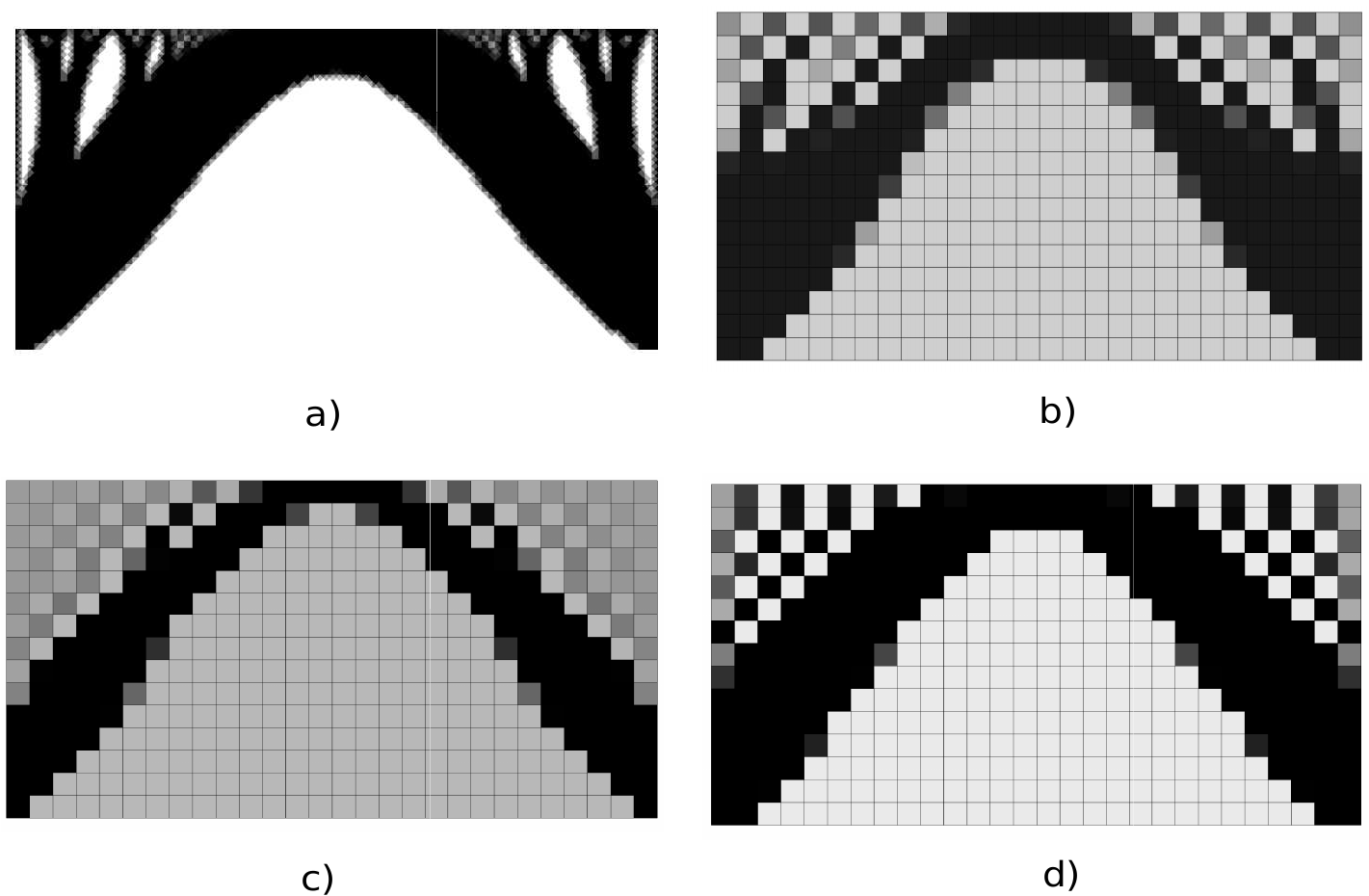}
\caption{Density distribution for the arch benchmark. a) TopOpt b) Lattice cell c) Circular inclusion d) Square inclusion. For this benchmark, the gray areas are largely reduced.}
\label{ArchDensities}
\end{figure}

Figure \ref{archTopolgies}b shows the topology obtained using the lattice cell parameterization for this benchmark. As in the cantilever beam case, we observe a global load path similar to that obtained with classical TopOpt (Figure \ref{archTopolgies}a), characterized by diagonal load-transfer patterns across the domain and vertical members descending from the top surface. However, gray regions reduce in size, with a clear tendency towards a black-and-white layout (see Figure \ref{ArchDensities}b).

Similarly, the structure resulting from the circular inclusion (Figure \ref{archTopolgies}c) exhibits a dominant diagonal load path together with weaker secondary diagonals above the main one. In this case, the topology tends toward a more black-and-white distribution as for the lattice-cell solution.

The square inclusion reproduces the main diagonal load path and also generates two additional smaller diagonals above the primary one (Figure \ref{archTopolgies}d). We also observe that the branching near the top boundary is consistent with the topology obtained using classical TopOpt.

\begin{figure}[hbt!] 
\centering
\includegraphics[width=0.88\columnwidth]{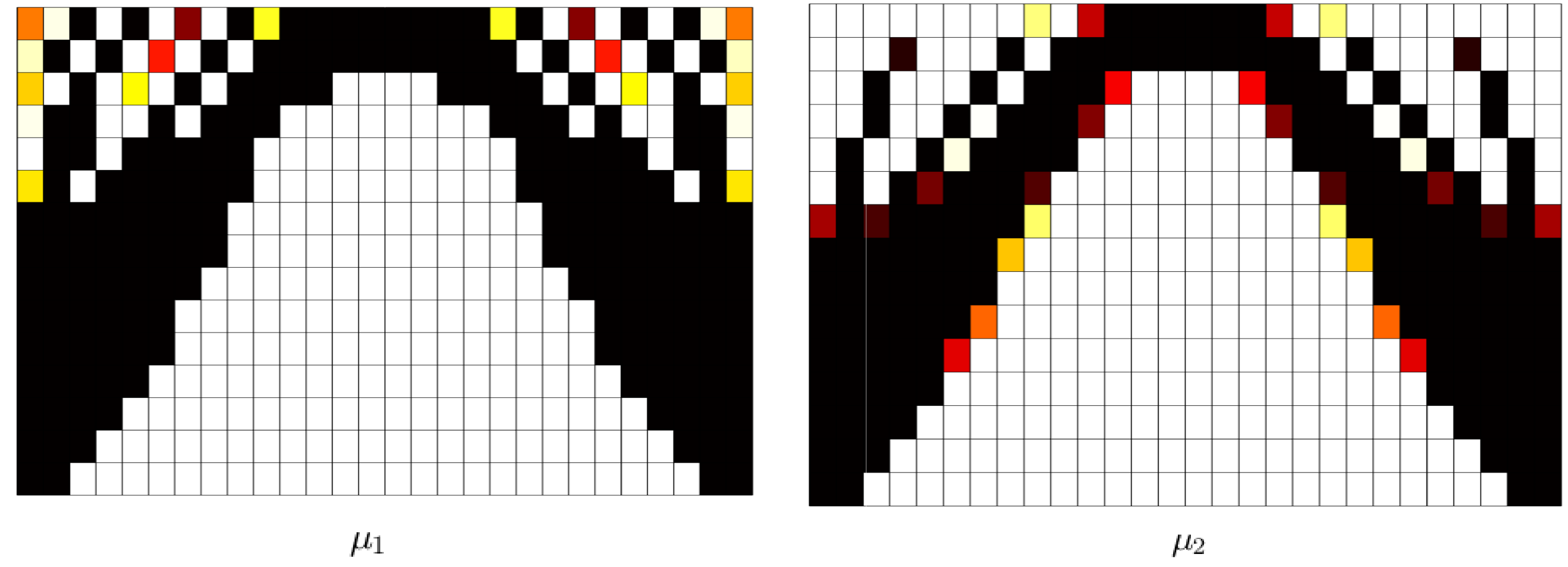}
\caption{Parameter distribution for the lattice cell. Black maximum value, white minimum value. When the bars are not aligned with the principal stresses, the results show a tendency towards the extreme values, removing intermediate ones and yielding more flexible structures.}
\label{ArchParameters}
\end{figure}

\begin{figure}[htbp]
    \centering
    \begin{subfigure}{0.58\textwidth}
        \includegraphics[width=\textwidth]{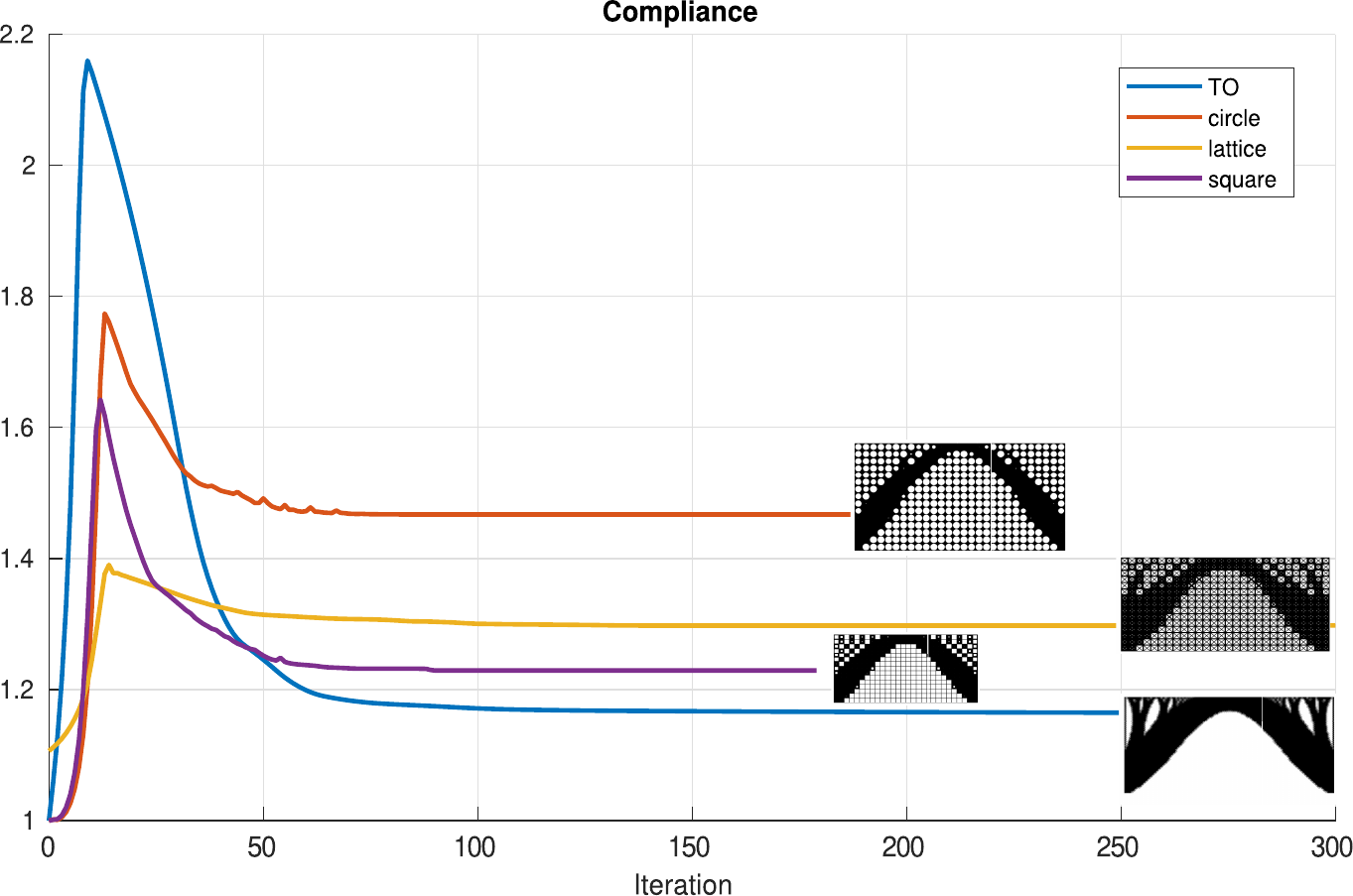}
        \caption{Convergence plots of the compliance for the arch benchmark. Results highlight the superior stiffness obtained via classical TopOpt compared to the parameterized microstructures. Within the parametric methods, the square inclusion demonstrates the most effective material redistribution.}
        \label{archMonitoring}
    \end{subfigure}
    \quad
    \begin{subfigure}{0.72\textwidth}
        \includegraphics[width=\textwidth]{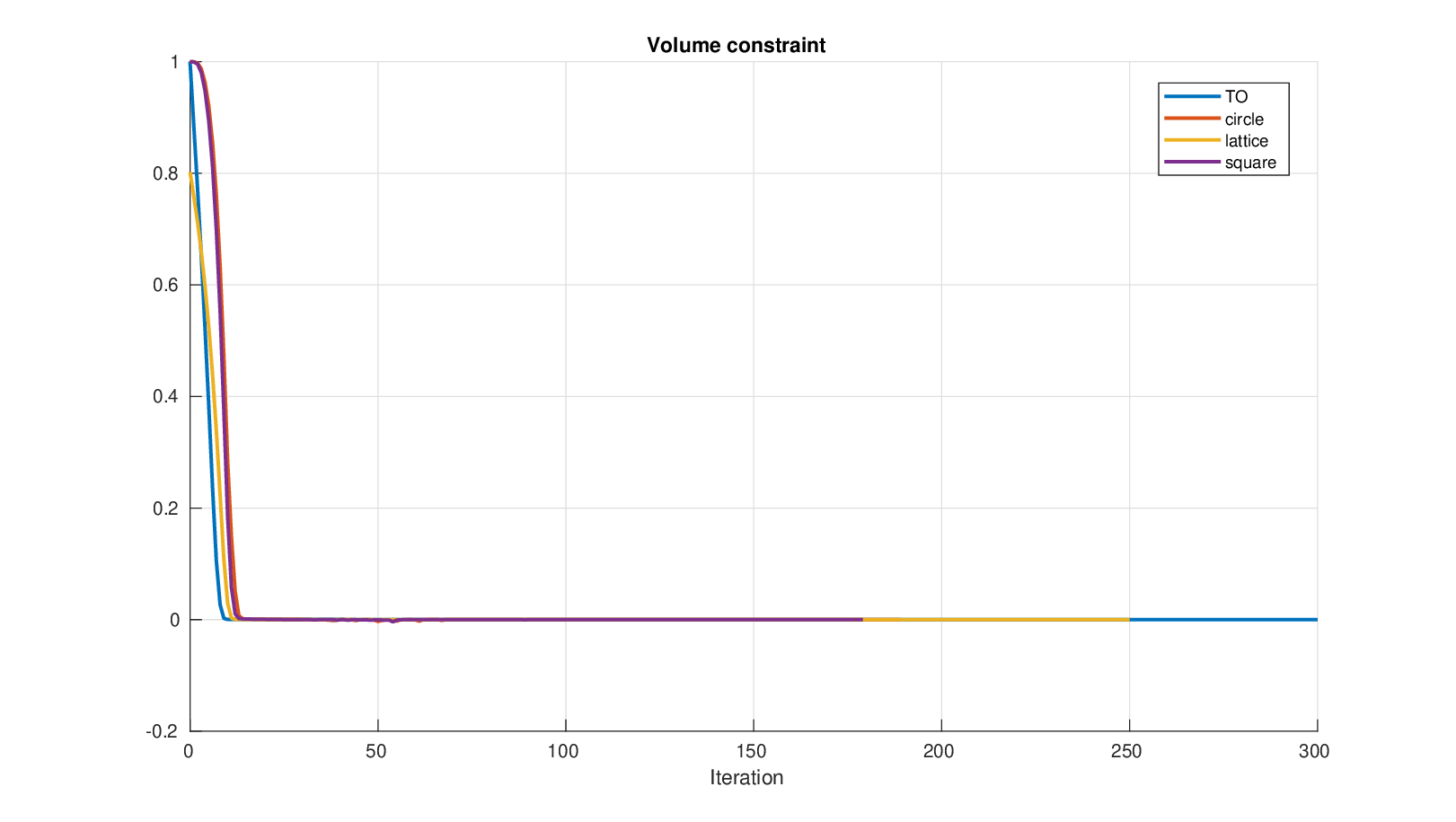}
        \caption{Relative volume constraint evolution as a function of the iteration.}
        \label{volumeArch}
    \end{subfigure}
    
    \caption{ Arch results. (a) Compliance evolution for TopOpt and the 3 different unit cells and (b) Volume evolution as a function of the iterations}
    \label{Arch results}
\end{figure}

Figure \ref{archMonitoring} shows the evolution of compliance for the four cases considered and Figure \ref{volumeArch} the evolution of the volume constraint. In contrast to the cantilever example, the stiffest structure is obtained using the classical TopOpt algorithm. Among the parametric approaches, the square inclusion yields the most competitive design.

This behavior can be explained by how the load is transferred in the arch benchmark.  The structural response in this domain is mainly governed by diagonal compressive load paths. Classical TopOpt can efficiently redistribute material along these continuous diagonal force-transfer paths, resulting in the stiffest overall structure.

The difference in performance between classical TopOpt and the lattice cell arises from the orientation of these stress trajectories. The diagonal compression paths in the arch benchmark do not align with the anisotropy introduced by the microstructure yielding more flexible structures. That gray areas reduce in size is also explained by this fact, since the most optimal way to stiffen the structure is by maximizing the thickness of both bar and frame (see Figure \ref{ArchParameters}). 

Regarding the square inclusion, it can redistribute mass more effectively where it is most advantageous; hence, for compression-dominated loads, it results in stiffer structures than the lattice cell. However, its ability to distribute mass is still less efficient than classical TopOpt, yielding higher values of compliance. 

\clearpage

\subsection{MBB}
The last example is geometrically different, since the aspect ratio of the domain is set to 4. This translates into a geometry formed by $40 \times10$ domains and 170.000 nodes. In the TopOpt case, the number of nodes in the discretization is 19750. The volume target is $50\%$.

\begin{figure}[hbt!] 
\centering
\includegraphics[width=0.65\columnwidth]{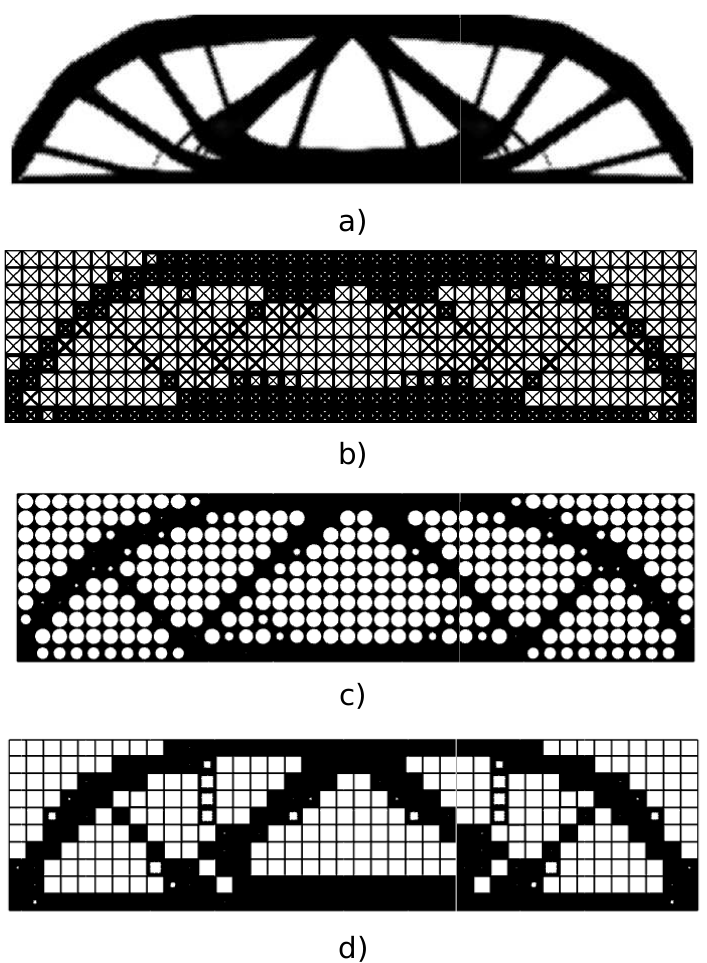}
\caption{MBB design obtained with a) TopOpt b) lattice unit cell c) circular and d) square inclusion. Parametric optimization obtains qualitatively similar structural pattern to classical TopOpt. Similarly to the cantilever beam, regions with intermediate densities appear again, especially for the lattice unit cell. This is a consequence of the cross-bars being aligned with the principal stresses which promotes the cross-bars to grow independently od the frame where advantageous.}
\label{MBBTopologies}
\end{figure}

Figure \ref{MBBTopologies} shows the final design for all 4 strategies. We observe a similar pattern between the TopOpt baseline and the parametric optimization, with an external arch to withstand bending connected by diagonal trusses.  

\begin{figure}[htbp]
    \centering
    \begin{subfigure}{0.58\textwidth}
        \includegraphics[width=\textwidth]{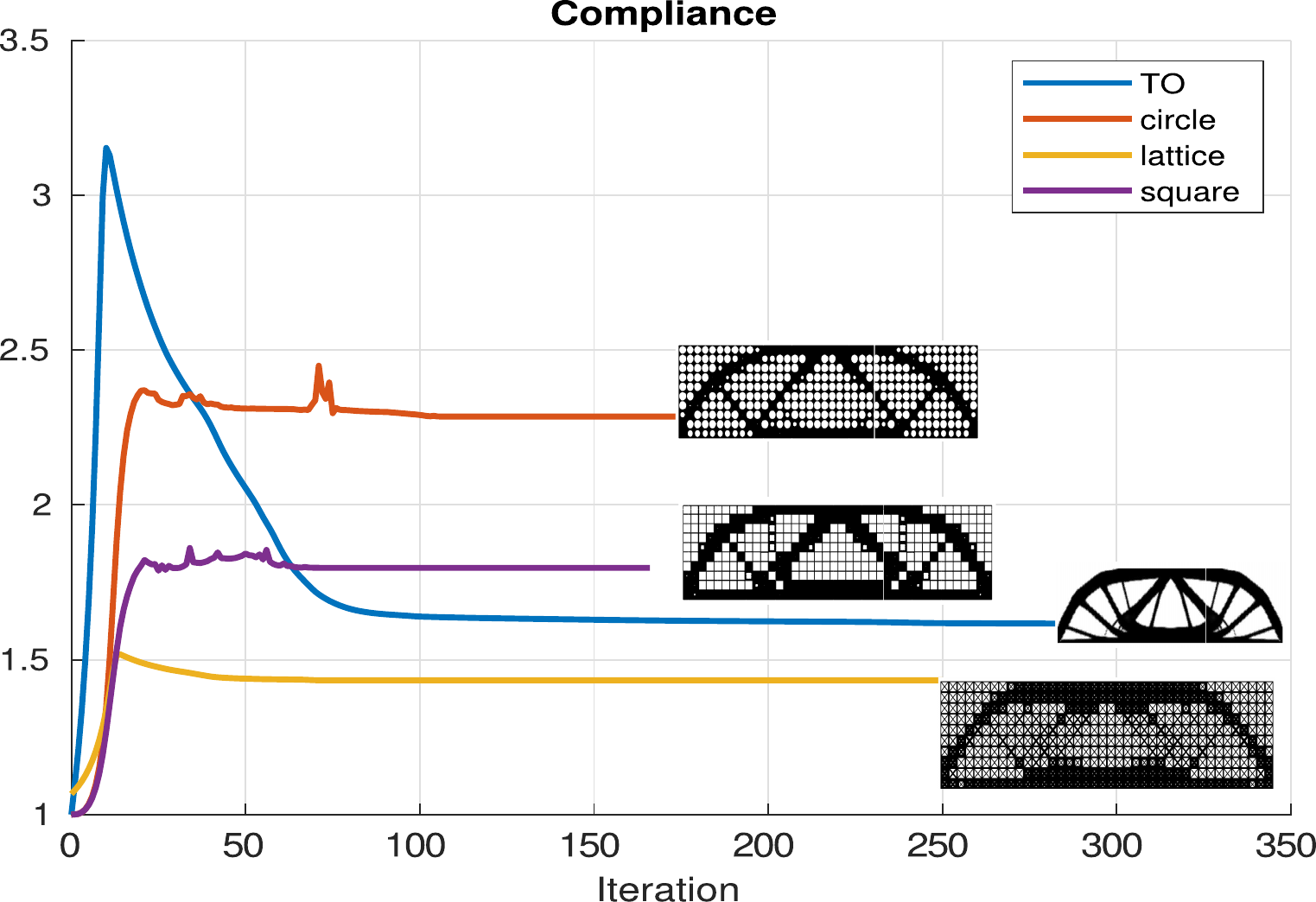}
        \caption{Convergence plots of the compliance objective function for the MBB benchmark. Results highlight the superior stiffness obtained with the lattice cell compared to classical TopOpt. Within the single parametric geometries, the square inclusion outperforms the circular inclusion. These results are inline with the cantilever benchmark.}
        \label{MBBMonitoring}
    \end{subfigure}
    \quad
    \begin{subfigure}{0.72\textwidth}
        \includegraphics[width=\textwidth]{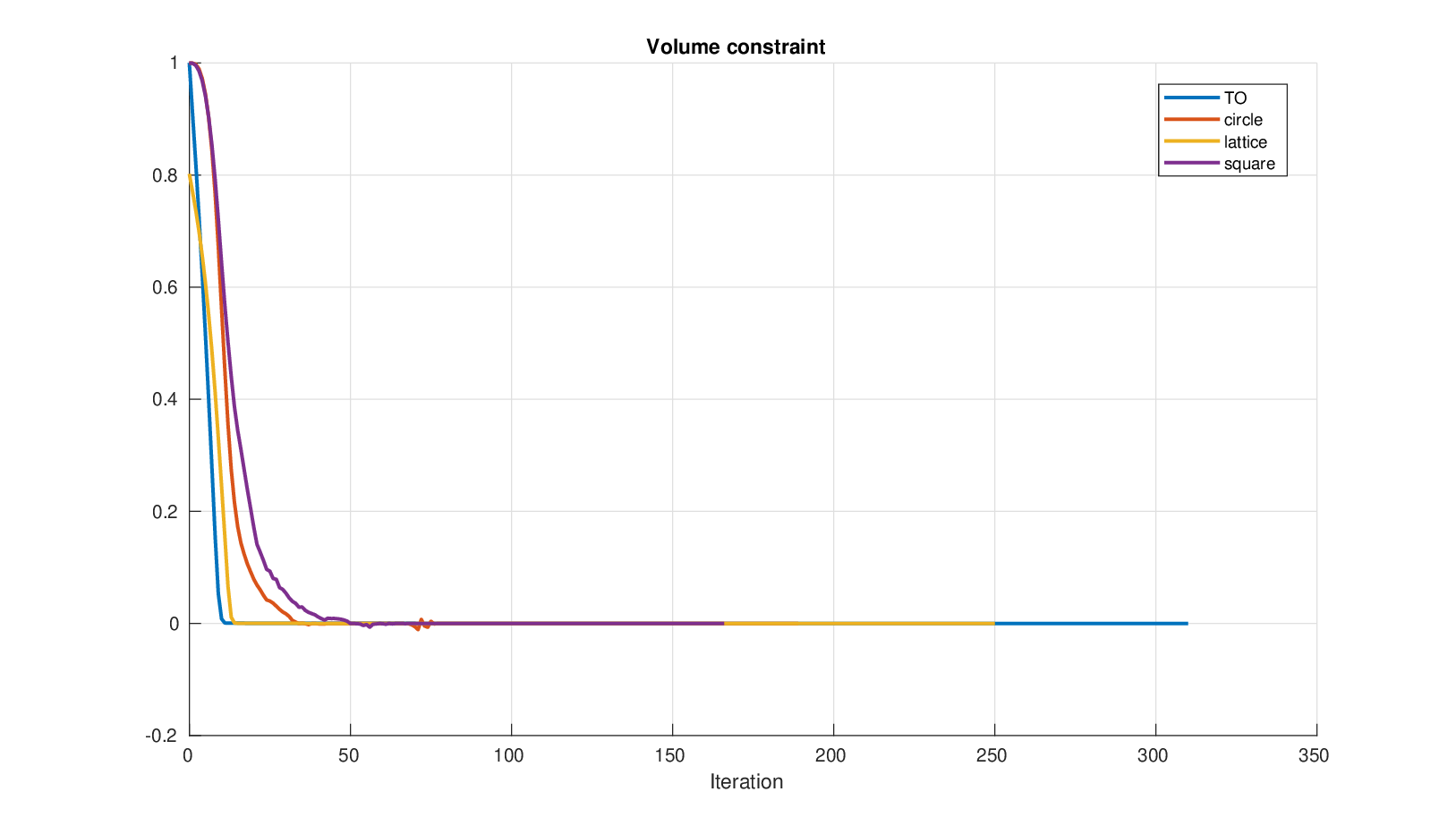}
        \caption{Relative volume constraint evolution as a function of the iteration.}
        \label{volumeMBB}
    \end{subfigure}
    
    \caption{ MBB results. (a) Compliance evolution for TopOpt and the 3 different unit cells and (b) Volume evolution as a function of the iterations}
    \label{MBB results}
\end{figure}

Figure \ref{MBBMonitoring} shows the evolution of compliance and the final optimized topologies for the MBB benchmark. Consistent with the cantilever results, the lattice parameterization achieves the lowest compliance, outperforming the classical TopOpt approach. Among the remaining parametric designs, the square inclusion performs better than the circular one, though both fall short of classical TopOpt.

\begin{figure}[hbt!] 
\centering
\includegraphics[width=0.88\columnwidth]{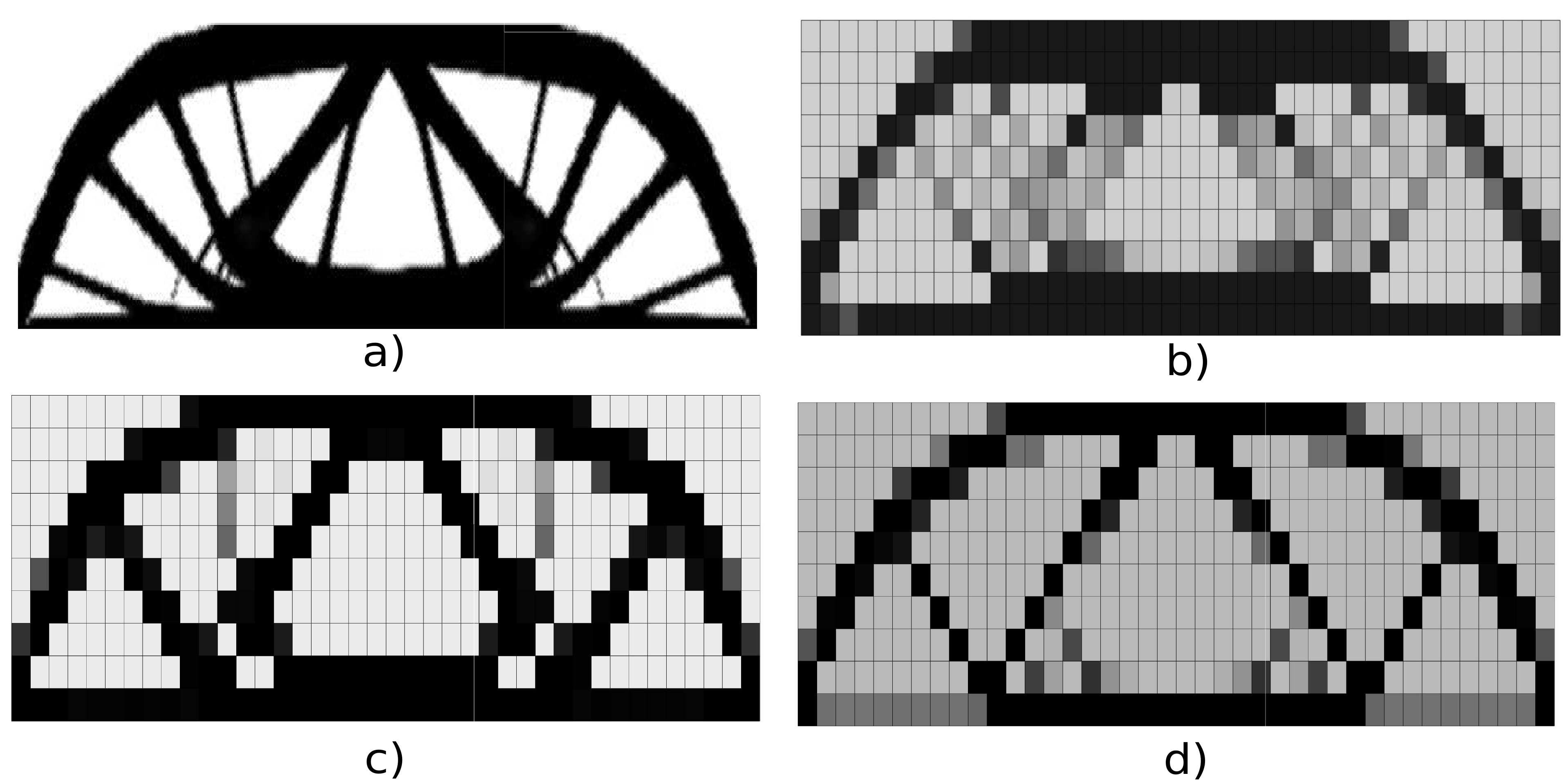}
\caption{Density distribution for the MBB benchmark. a) TopOpt b) Lattice cell c) Square inclusion d) Circular inclusion. For this benchmark, the gray areas we observed in the cantilever example appear again. This is due to the shear bands produced by the loading, which promotes the growth of the cross-bars (see Figure \ref{MBBParameters}).}
\label{MBBDensities}
\end{figure}

\begin{figure}[hbt!] 
\centering
\includegraphics[width=0.88\columnwidth]{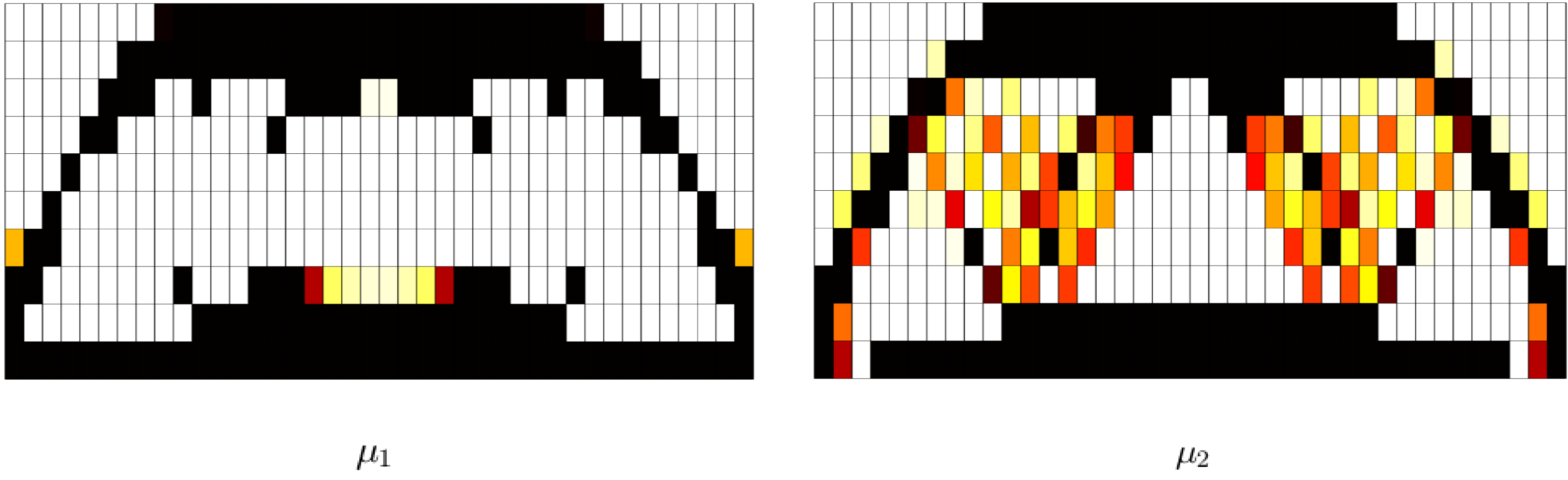}
\caption{Parameter distribution for the lattice cell. Black maximum value, white minimum value. As in the cantilever beam example, this loading configuration produces shear bands. To withstand the stresses, cross bars are very effective and thus the optimizer adds intermediate values regions for the parameter $\mu_2$.}
\label{MBBParameters}
\end{figure}

The better performance of the lattice cell in this example, can be attributed (again) to how well its anisotropy aligns with the stress trajectories. The structural response of an MBB beam is mainly governed by horizontal tension and compression along the outer boundaries (global bending) and diagonal shear paths throughout the web. The lattice cell, which consists of both an outer frame and internal diagonal bars, is able to transmit this stress configuration very effectively. Thus, the "gray" areas become highly effective in stiffening the structure, as shown in Figure \ref{MBBDensities}b, and result directly from the parameter distribution (Figure \ref{MBBParameters}).

Classical TopOpt yields the second-stiffest structure. Because it is not constrained by a predefined shape, it distributes solid material along the optimal diagonal shear paths. However, it lacks the mechanical advantage of the intermediate-density regions exploited by the lattice.

Finally, the performance gap between the square and circular inclusions highlights the importance of geometric flexibility. Although the square inclusion layout is well-suited for the horizontal bending stresses, it lacks the diagonal members necessary to efficiently manage shear, leading to a higher compliance than the lattice. The circular inclusion, which has a lower ability to redistribute mass than the square inclusion and TopOpt, and lacks distinct directional bars, produces the most flexible structure.

\clearpage

\subsection{Computational time comparison}

Table \ref{tab:time comparison} summarizes the computational time detailing the offline and online costs across all three benchmarks: the cantilever, arch, and MBB beam. The data reveals a consistent and remarkable online speedup achieved by the parametric EIFEM-based optimization regardless of the structural domain.

The parametric models solve fine-scale systems (over 330,000 DoF) in a matter of seconds (ranging from 4.5 s to 14 s). In contrast, classical TopOpt requires several minutes to solve significantly coarser systems. Furthermore, the performance gap widens as the problem size increases. For example, in the MBB benchmark, the classical TopOpt solve time jumps to nearly 6 minutes for a 39,500 DoF system. Meanwhile, the lattice parameterization solves a highly detailed 330,000 DoF system in just 7.6 seconds.

While the lattice cell demands the highest offline computational cost (10 minutes) to construct the model compared to the single parameter inclusions (31 seconds), this is a fixed, one-time cost. Given the massive reduction in iterative online solve and the superior mechanical performance of the lattice demonstrated in the cantilever and MBB examples, this initial offline is highly justifiable.

\begin{table}[h]
\centering
\makebox[\textwidth][c]{%
\begin{tabular}{ll|ll|ll|ll|}
\cline{3-8}
                              &              & \multicolumn{2}{c|}{Cantilever} & \multicolumn{2}{c|}{Arch} & \multicolumn{2}{c|}{MBB} \\ \cline{2-8} 
\multicolumn{1}{l|}{}         &\begin{tabular}{c}
Offline
time
\end{tabular} & Dofs         & Time       & Dofs      & Time    & Dofs      & Time  \\ \hline
\multicolumn{1}{|l|}{Circle}  & 31 s         & 340000       & 4.5 s            & 340000    & 7 s           & 330000    & 6.1 s        \\
\multicolumn{1}{|l|}{Square}  & 31 s         & 340000       & 5.6 s            & 340000    & 8.4 s         & 330000    & 6.8 s        \\
\multicolumn{1}{|l|}{Lattice} & 10 min       & 340000       & 13 s             & 340000    & 14 s          & 330000    & 7.6 s        \\
\multicolumn{1}{|l|}{TopOpt}  & $\emptyset$  & 19600        & 2 min            & 19600     & 2 min 39 s    & 39500     & 5 min 45 s   \\ \hline
\end{tabular}}
\caption{Summary of degrees of freedom, offline model construction time, and online solve times for the considered benchmarks. The parametric methods consistently achieve online solve times below 15 seconds for fine-scale microstructural meshes (>330,000 DoFs). While the lattice parameterization has a 10-minute offline cost, it offers remarkable online speedups compared to classical TopOpt, particularly as the macroscopic domain size increases (e.g., the MBB beam).}
\label{tab:time comparison}
\end{table}

\subsection{3D Cantilever beam}

\begin{figure}[hbt!] 
\centering
\makebox[\columnwidth][c]{%
\includegraphics[width=0.9\columnwidth]{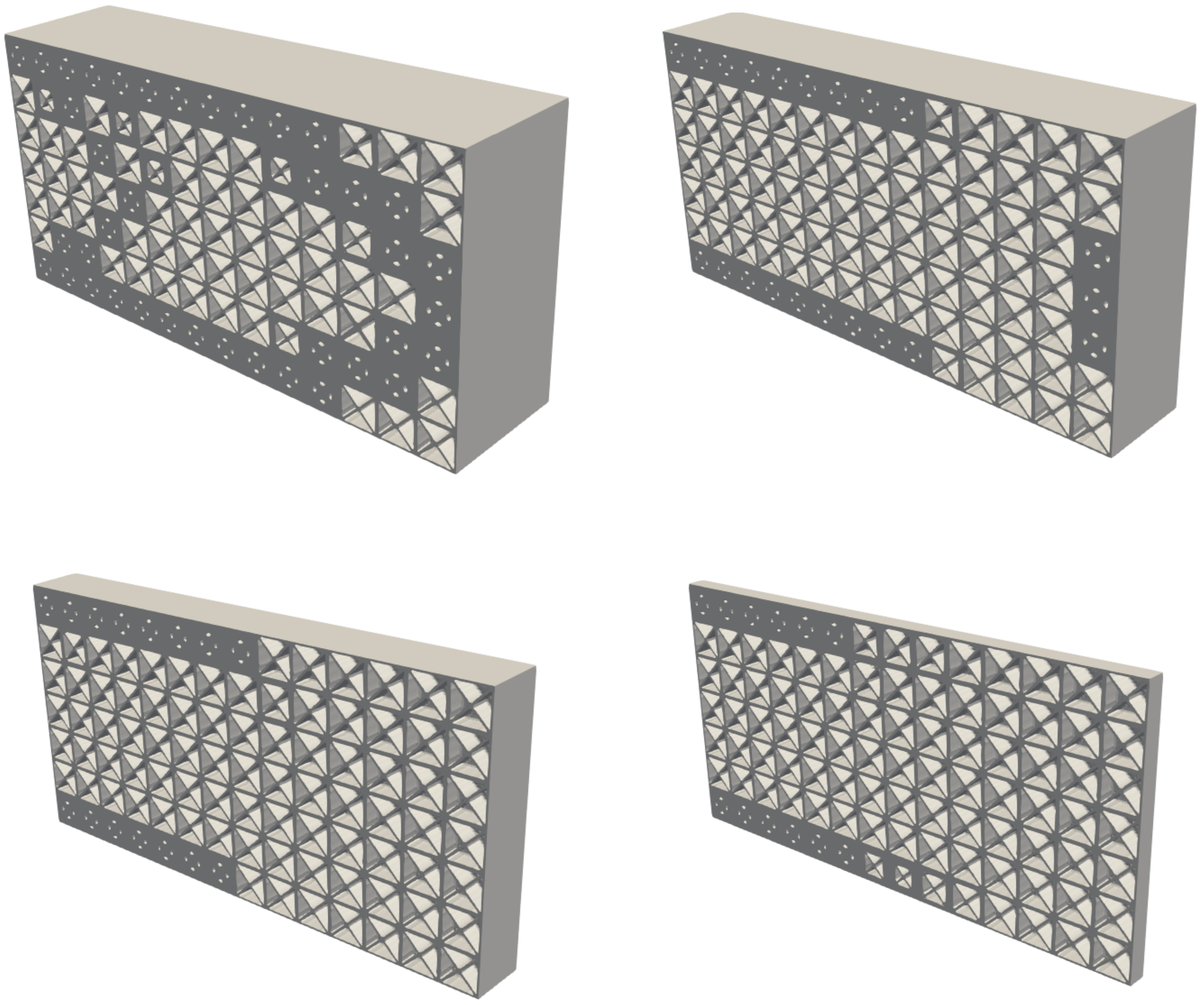}
}
\caption{Topology of the 3D cantilever beam. Each image represents a section of the geometry, from the middle plane to the outer skin. The location of the cuts are selected such that the middle section of the cell is visible.}
\label{3DLatticeCuts}
\end{figure}

\begin{figure}[hbt!] 
\centering
\makebox[\columnwidth][c]{%
\includegraphics[width=1\columnwidth]{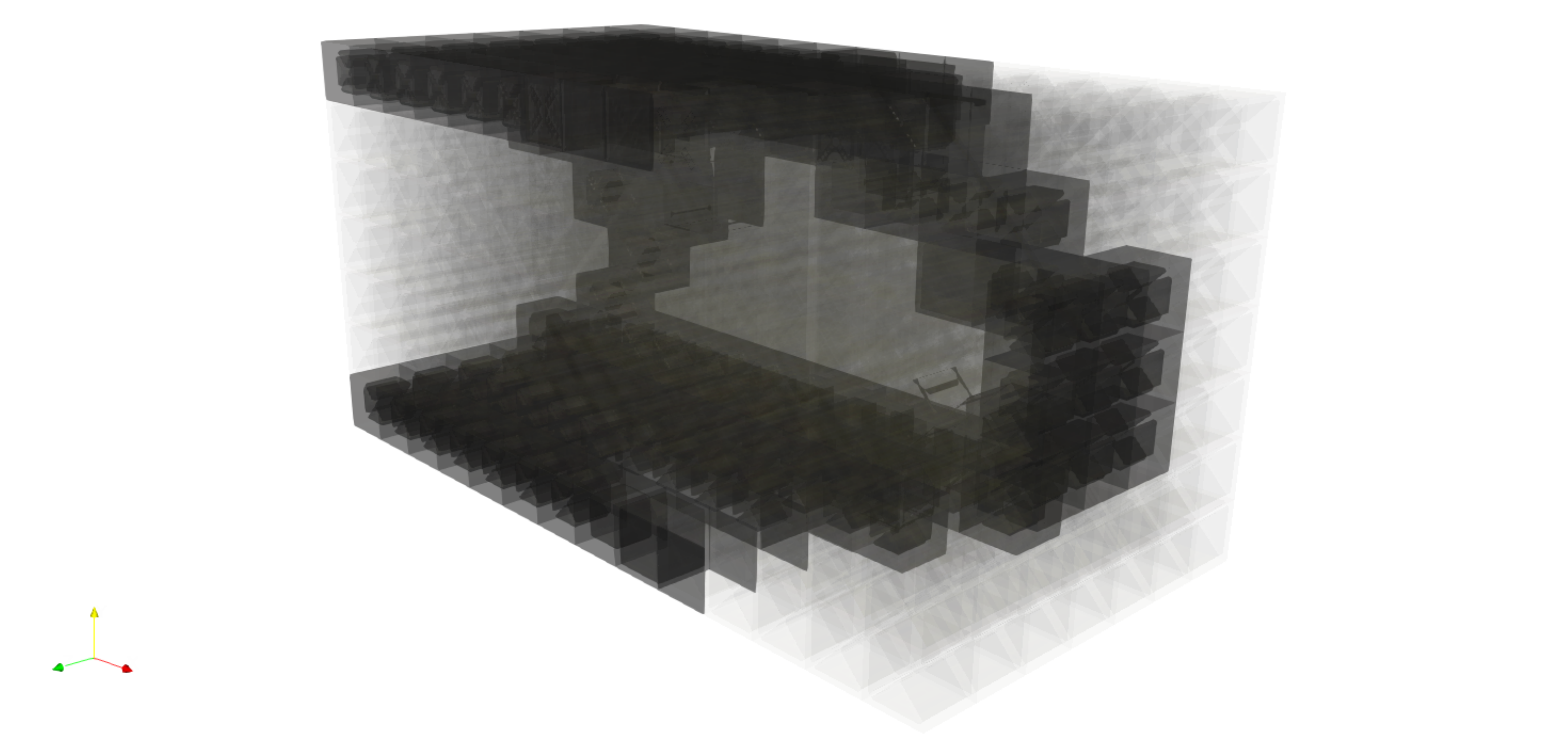}
}
\caption{Topology of the 3D cantilever beam. Transparency added just for visualization purposes. This image complement Figure \ref{3DLatticeCuts}.}
\label{3DLatticeTopology}
\end{figure}

\begin{figure}[hbt!]
    \centering
    \begin{subfigure}{0.8\textwidth}
       \makebox[\columnwidth][c]{%
 \includegraphics[width=1\columnwidth]{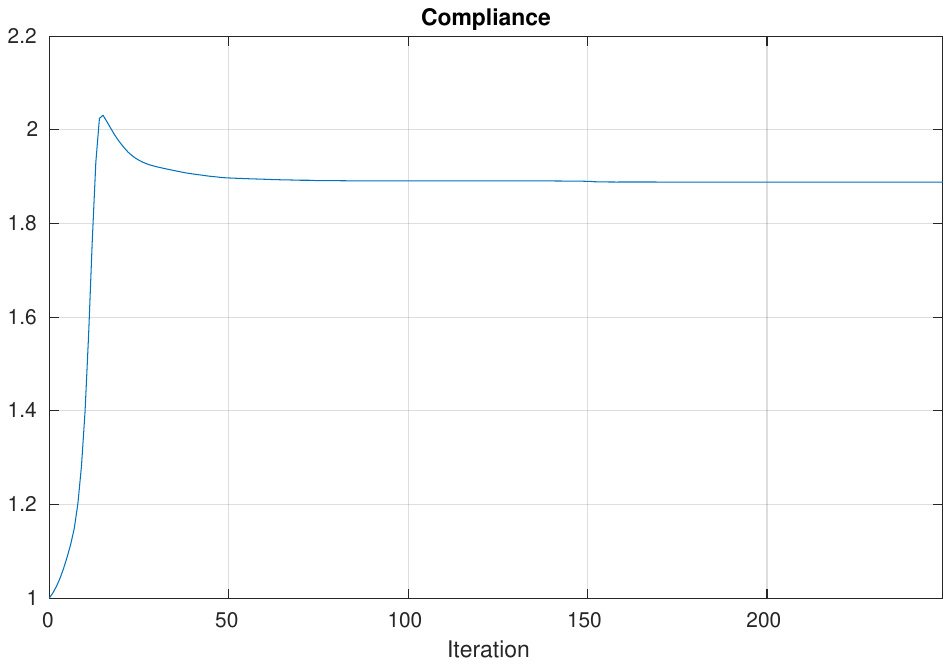}
 }
        \vspace{-145pt}
        \caption{Compliance evolution as a function of the iteration.}
        \label{3Dmonitoring}
    \end{subfigure}
    \begin{subfigure}{0.72\textwidth}
        \includegraphics[width=\textwidth]{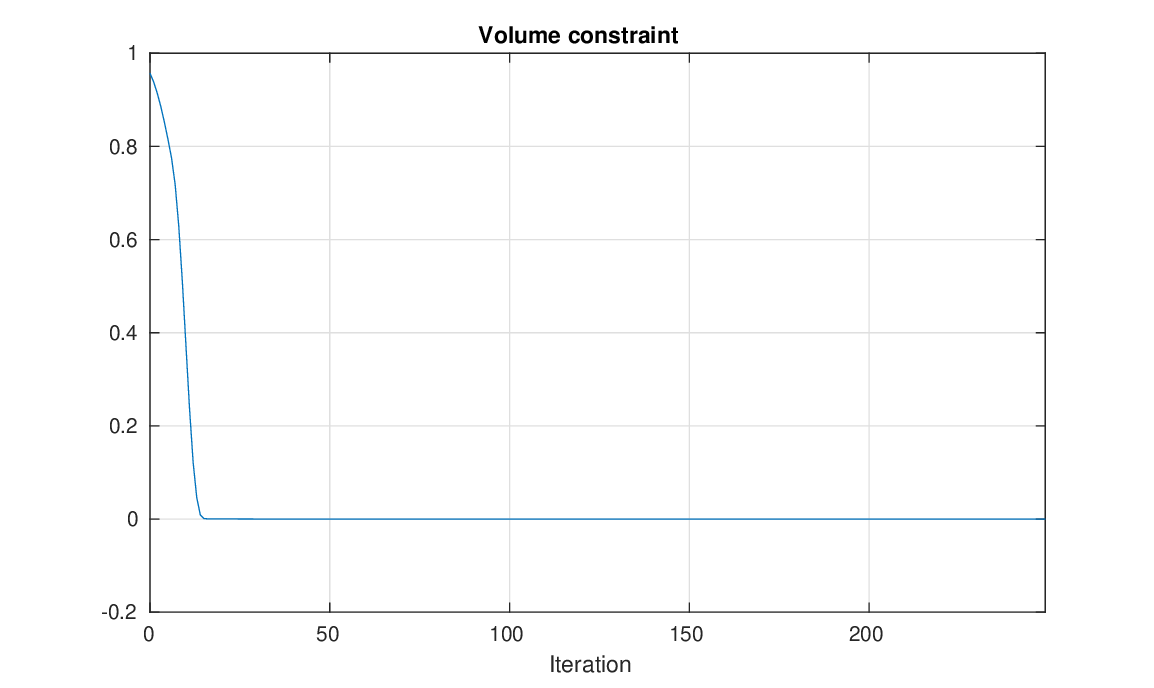}
        \caption{Relative volume constraint evolution as a function of the iteration.}
        \label{volume3D}
    \end{subfigure}
    
    \caption{ Monitoring of the 3D cantilever beam with the lattice cell. Optimizer MMA.}
    \label{3D results}
\end{figure}

To demonstrate the scalability and practical applicability of the proposed methodology, the parametric ROM approach was extended to a 3D cantilever benchmark. Figure \ref{3DLatticeCuts} and Figure \ref{3DLatticeTopology} illustrate the final optimized topology using the 3D lattice parameterization, displaying several cross-sections from the middle plane to the outer skin. The evolution of the compliance is shown in Figure \ref{3Dmonitoring} while that of the volume in Figure \ref{volume3D}.

The resulting geometry highlights the mechanical efficiency of the lattice microstructure in three dimensions. The optimizer naturally distributes mass to form a sandwich-like composite structure. The outer skins (top and bottom regions) have higher concentrations of solid material to resist the maximum tensile and compressive bending stresses. In contrast, the internal core is mainly populated by intermediate-density lattice cells. These internal layouts efficiently transfer 3D shear forces across the domain while drastically reducing the overall weight.

Consistent with the 2D findings, this 3D result highlights the structural advantage of exploiting "gray" regions rather than penalizing them. By physically realizing these intermediate densities as an oriented 3D lattice, the structure benefits from geometry aligned with the spatial stress trajectories. Furthermore, while full-scale 3D classical topology optimization typically demands massive computational resources, the parametric ROM framework allows this complex 3D structure, with 1.8 million DoFs, to be solved in less than 15 seconds during the online optimization phase.

\clearpage

\clearpage

\section{Conclusions} \label{Conclusions}
This article has presented and evaluated a parametric reduced-order modeling strategy for structural optimization within the EIFEM framework. The resulting algorithm can be regarded as a preconditioner for the optimization problem, where both the design variable and the state equation are approximated with coarse spaces.

The central ingredient is the efficient approximation of the parameter-dependent interscale operator and the associated coarse stiffness matrix, which must be repeatedly evaluated during the optimization loop. To this end, the Discrete Empirical Interpolation Method (DEIM) is used to express these high-dimensional operators in terms of a reduced set of representative entries selected from offline snapshots. These selected components are then interpolated with respect to the geometric parameters, and the full operators are reconstructed from this reduced information. In this way, the dependence of the EIFEM operators on the design variables can be captured at a significantly reduced cost, without assembling the full high-fidelity matrices at each iteration.

From a numerical standpoint, the results show that this strategy yields accurate operator reconstructions, even when using relatively sparse sampling in the parameter space, and does not degrade the performance of the preconditioning strategy introduced in \cite{RUBIO2025}. In addition, the approach can be integrated into existing finite element codes with minimal modifications, as it relies primarily on the evaluation of parametric operators through function handles.

Overall, this contribution demonstrates how the EIFEM framework can be extended to parametric settings in a consistent offline/online fashion, enabling efficient optimization while preserving the modular and data-driven nature of the underlying methodology.


Mechanically, the benchmark results highlight the structural value of intermediate-density ("gray") regions, which are typically penalized and removed in classical TopOpt. When these gray areas are physically realized as a lattice, they can produce highly competitive designs. The overall stiffness of the optimized structure is strongly dependent on how well the microstructure aligns with the dominant stress trajectories. These findings align with the principles of dehomogenization, where macroscopic intermediate densities are mapped to optimally oriented microstructures. Conversely, isotropic microstructures consistently perform worse than classical TopOpt. They lack both directional stiffness and geometric flexibility to effectively redistribute mass, essentially acting as a coarse-resolution TopOpt that is too restricted to form optimal load paths.

Finally, while the ROM-based parametrization restricts the design space compared to the full topological freedom of density-based methods and requires a non-negligible offline computational investment, these limitations are compensated by significant gains in the online stage. In particular, by combining EIFEM with DEIM, the computational cost of the optimization loop is reduced by one to two orders of magnitude compared to standard approaches, as the repeated assembly and solution of high-fidelity models is avoided.

Moreover, unlike density-based topology optimization, the proposed framework directly yields manufacturable designs described by explicit geometric parameters, eliminating the need for additional post-processing steps. Overall, for applications requiring rapid design iterations, the proposed approach provides a scalable and computationally efficient alternative, effectively shifting the computational effort to the offline stage while enabling fast evaluations during optimization.

\section*{Declaration of Generative AI and AI-assisted Technologies in the Writing Process}
During the preparation of this manuscript, the authors used ChatGPT and GEMINI to assist with improving the grammar, clarity, and style of the text. The AI system was not used to generate scientific results, perform data analysis, develop the methodology, or draw conclusions. All scientific content, interpretations, and final editorial decisions were reviewed and verified by the authors, who take full responsibility for the content of the manuscript.

\section*{Acknowledgement}
Part of this research was conducted during a research stay at École polytechnique fédérale de Lausanne (EPFL). R.Rubio sincerely acknowledges the support and stimulating research environment provided by the institution.

R.Rubio acknowledges the support of Beques Santander – Ayudas predoctorales 2025.

R.Rubio and A.Ferrer gratefully acknowledge the support of the following projects: "Optimización topológica con compuestos y materiales flexibles" (PID2023-153213NA-I00), "Disseny computacional òptim de bateries amb materials compostos sostenibles per a una mobilitat més eficient" (ACE122/24/000030) and "Design of flexible materials, structures via shape and topology optimization" (PCI2024-155060-2).

J.A. Hernández acknowledges the support of Grant PID2024-158878OB-C21 funded by MICIU/AEI /10.13039/501100011033 and by ERDF/EU.

P.Antolin acknowledges the financial support of the Swiss National Science Foundation through the project FLAS\emph{h} with no. 200021\_214987

\bibliography{references}

@article{Park2000,
title = {A variational principle for the formulation of partitioned structural systems},
journal = {Numerical Methods in Engineering},
volume = {47},
pages = {395-418},
year = {2000},
issn = {0045-7825},
doi = {10.1002/(SICI)1097-0207(20000110/30)47:1/3<395::AID-NME777>3.0.CO;2-9},
author = {Kwang-Chung Park and Carlos Alberto Felippa},
}

@article{DVORAK2024,
title = {On the automatic construction of interface coupling operators for non-matching meshes by optimization methods},
journal = {Computer Methods in Applied Mechanics and Engineering},
volume = {432},
pages = {117336},
year = {2024},
issn = {0045-7825},
doi = {10.1016/j.cma.2024.117336},
author = {R. Dvořák and J A González}
}

@book{TosselliWidlund2006,
    title={Domain Decomposition Methods - Algorithms and Theory },
    author = { Andrea Toselli and Olof B. Widlund},
    year={2006},
    publisher={Springer}
}

@book{Doolean2015,
    title={An Introduction to Domain Decomposition Methods: algorithms, theory and parallel implementation},
    author = {Victorita Doolean and Pierre Jolivet and Frédéric Nataf},
    year={2015},
    publisher={SIAM}
}

@book{Trottenberg2000,
  title={Multigrid},
  author={Trottenberg, U. and Oosterlee, C.W. and Schuller, A.},
  isbn={9780080479569},
  year={2000},
  publisher={Academic Press}
}

@book{Nocedal2006,
  title={Numerical optimization},
  author={Jorge Nocedal and J S Smith},
  year={2006},
  publisher={Springer},
}

@book{Olshanskii2014,
  title={Iterative Methods for Linear Systems: Theory and Applications},
  author={Maxim A. Olshanskii and Eugene E. Tyrtyshnikov},
  year={2014},
  publisher={SIAM}
}

@article{Farhat2001FETIDP,
  author = {Charbel Farhat and Michel Lesoinne and Patrick Le Tallec and Kendall Pierson and Daniel Rixen},
  title = {{FETI-DP: A Dual–Primal Unified FETI Method—Part I: A Faster Alternative to the Two-Level FETI Method}},
  journal = {International Journal for Numerical Methods in Engineering},
  volume = {50},
  number = {7},
  pages = {1523--1544},
  year = {2001},
  doi = {10.1002/nme.116},
}

@article{mandel2008multispace,
  title={Multispace and multilevel BDDC},
  author={Mandel, Jan and Sousedík, Bedřich and Clark R. Dohrmann },
  journal={Computing},
  volume={83},
  number={2},
  pages={55-85},
  year={2008},
  publisher={Springer},
  doi = {10.1007/s00607-008-0014-7}
}

@book{Briggs2000,
author = {Briggs, William L. and Henson, Van Emden and McCormick, Steve F.},
title = {A Multigrid Tutorial, Second Edition},
publisher = {Society for Industrial and Applied Mathematics},
year = {2000},
doi = {10.1137/1.9780898719505},
edition   = {Second},
}

@article{Berkooz1993,
   author = "Berkooz, G and Holmes, P and Lumley, J L",
   title = "The Proper Orthogonal Decomposition in the Analysis of Turbulent Flows", 
   journal= "Annual Review of Fluid Mechanics",
   year = "1993",
   volume = "25",
   number = "Volume 25, 1993",
   pages = "539-575",
   doi = "10.1146/annurev.fl.25.010193.002543"
  }

@article{RUBIO2025,
title = {Preconditioning iterative solvers via the {E}mpirical {I}nterscale {F}inite {E}lement {M}ethod ({EIFEM})},
journal = {Computer Methods in Applied Mechanics and Engineering},
volume = {446},
pages = {118257},
year = {2025},
issn = {0045-7825},
doi = {10.1016/j.cma.2025.118257},
author = {R. Rubio and A. Ferrer and J.A. Hernández}
}

@article{Chaturantabut2010,
author = {Chaturantabut, Saifon and Sorensen, Danny C.},
title = {Nonlinear Model Reduction via Discrete Empirical Interpolation},
journal = {SIAM Journal on Scientific Computing},
volume = {32},
number = {5},
pages = {2737-2764},
year = {2010},
doi = {10.1137/090766498},
}

@book{quarteroni2016reduced,
  title     = {Reduced Basis Methods for Partial Differential Equations: An Introduction},
  author    = {Quarteroni, Alfio and Manzoni, Andrea and Negri, Federico},
  year      = {2016},
  publisher = {Springer Cham},
  series    = {UNITEXT},
  isbn      = {978-3-319-15430-5, 978-3-319-15431-2},
  doi       = {10.1007/978-3-319-15431-2}
}

@book{Canuto2006,
  author    = {Canuto, Claudio and Hussaini, M. Yousuff and Quarteroni, Alfio and Zang, Thomas A.},
  title     = {Spectral Methods: Fundamentals in Single Domains},
  publisher = {Springer},
  year      = {2006},
  series    = {Scientific Computation},
  address   = {Berlin, Heidelberg},
  isbn      = {978-3-540-30726-6}
}

@book{Hesthaven_2007, 
title={Spectral Methods for Time-Dependent Problems},
publisher={Cambridge University Press},
author={Hesthaven, Jan S. and Gottlieb, Sigal and Gottlieb, David},
year={2007}
}

@book{BendsoeSigmund2003,
  author    = {Bend{\o}e, Martin Philip and Sigmund, Ole},
  title     = {Topology Optimization: Theory, Methods, and Applications},
  year      = {2003},
  publisher = {Springer Science \& Business Media},
  address   = {Berlin, Heidelberg},
  isbn      = {978-3-540-42992-1}
}

@article{Allaire2004,
  author  = {Allaire, Gr{\'e}goire and Jouve, Fran{\c{c}}ois and Toader, Anca-Maria},
  title   = {Structural optimization using sensitivity analysis and a level-set method},
  journal = {Journal of Computational Physics},
  volume  = {194},
  number  = {1},
  pages   = {363--393},
  year    = {2004},
  doi     = {10.1016/j.jcp.2003.09.032}
}

@article{Svanberg1987MMA,
  author  = {Svanberg, Krister},
  title   = {The method of moving asymptotes---a new method for structural optimization},
  journal = {International Journal for Numerical Methods in Engineering},
  year    = {1987},
  volume  = {24},
  number  = {2},
  pages   = {359--373},
  doi     = {10.1002/nme.1620240207}
}

@article{Feppon2020,
  author  = {Feppon, Fr{\'e}d{\'e}ric and Allaire, Gr{\'e}goire and Dapogny, Charles},
  title   = {Null space gradient flows for constrained optimization with applications to shape optimization},
  journal = {ESAIM: Control, Optimisation and Calculus of Variations},
  year    = {2020},
  volume  = {26},
  pages   = {90},
  doi     = {10.1051/cocv/2019026}
}

@article{PfluegerP2010,
  author       = {Dirk Pflüger and Benjamin Peherstorfer and Hans-Joachim Bungartz},
  title        = {Spatially Adaptive Sparse Grids for High-Dimensional Data-Driven Problems},
  journal      = {Journal of Complexity},
  volume       = {26},
  number       = {5},
  pages        = {508--522},
  year         = {2010},
  issn         = {0885-064X},
  doi          = {10.1016/j.jco.2010.04.001}
}

@article{Sigmund2013,
  author  = {Sigmund, Ole and Maute, Kurt},
  title   = {Topology optimization approaches: A comparative review},
  journal = {Structural and Multidisciplinary Optimization},
  year    = {2013},
  volume  = {48},
  number  = {6},
  pages   = {1031--1055},
  month   = {August},
  doi     = {10.1007/s00158-013-0978-6},
}

@book{Allaire2007,
  author    = {Gr{\'e}goire Allaire},
  title     = {Conception optimale de structures},
  series    = {Math{\'e}matiques et Applications},
  volume    = {58},
  publisher = {Springer-Verlag},
  address   = {Berlin},
  year      = {2007},
  isbn      = {978-3-540-69491-6}
}

@article{Lazarov2016,
  author    = {Lazarov, Boyan S. and Wang, Fengwen and Sigmund, Ole},
  title     = {Length scale and manufacturability in density-based topology optimization},
  journal   = {Archive of Applied Mechanics},
  year      = {2016},
  volume    = {86},
  number    = {},
  pages     = {189--218},
  doi       = {10.1007/s00419-015-1106-4}
}

@article{Amstutz2022,
  author    = {Amstutz, Samuel and Dapogny, Charles and Ferrer, Alex},
  title     = {A consistent approximation of the total perimeter functional for topology optimization algorithms},
  journal   = {ESAIM: Control, Optimisation and Calculus of Variations},
  year      = {2022},
  volume    = {28},
  pages     = {18},
  doi       = {10.1051/cocv/2022006}
}

@article{Allaire2017,
  author    = {Allaire, Gr{\'e}goire and Dapogny, Charles and Estevez, Romain and Faure, Adrien and Michailidis, Georgios},
  title     = {Structural optimization under overhang constraints imposed by additive manufacturing technologies},
  journal   = {Journal of Computational Physics},
  year      = {2017},
  volume    = {351},
  pages     = {295--328},
  doi       = {10.1016/j.jcp.2017.09.041}
}

@article{Torres2025,
author = {Torres, Jose and Otero, Fermin and Ferrer, Alex},
title = {Global Length and Overhang Control for Level Set and Density Approaches via Perimeter Minimization},
journal = {International Journal for Numerical Methods in Engineering},
volume = {126},
number = {2},
pages = {e7662},
doi = {10.1002/nme.7662},
year = {2025}
}

@article{Bayat2023,
  author    = {Bayat, M. and Zinovieva, O. and Ferrari, F. and Ayas, C. and Langelaar, M. and Spangenberg, J. and Salajeghe, R. and Poulios, K. and Mohanty, S. and Sigmund, Ole and others},
  title     = {Holistic computational design within additive manufacturing through topology optimization combined with multiphysics multi-scale materials and process modelling},
  journal   = {Progress in Materials Science},
  year      = {2023},
  pages     = {101129},
  doi       = {10.1016/j.pmatsci.2023.101129}
}

@incollection{Buhr2021,
title = {Localized model reduction for parameterized problems},
booktitle = {Volume 2: Snapshot-Based Methods and Algorithms},
author = {Buhr, Andreas and Iapichino, Laura and Ohlberger, Mario and Rave, Stephan and Schindler, Felix and Smetana, Kathrin},
editor = {Peter Benner and Stefano Grivet-Talocia and Alfio Quarteroni and Gianluigi Rozza and Wil Schilders and Luís Miguel Silveira},
publisher = {De Gruyter},
pages = {245--306},
doi = {10.1515/9783110671490-006},
isbn = {9783110671490},
year = {2021},
}

@phdthesis{Diercks2025,
title = "Multiscale modeling of mechanical structures via localized model reduction",
author = "Philipp Diercks",
note = "Proefschrift.",
year = "2025",
language = "English",
isbn = "978-90-386-6389-0",
publisher = "Eindhoven University of Technology",
type = "Phd Thesis",
school = "Mathematics and Computer Science",
}

@article{NEZDYUR2026,
title = {Parametric reduced order models for graded lattice structures},
journal = {Finite Elements in Analysis and Design},
volume = {255},
pages = {104521},
year = {2026},
issn = {0168-874X},
doi = {10.1016/j.finel.2026.104521},
author = {Max Nezdyur and Lynn Munday and Wilkins Aquino},
}

@article{Brezzi2005,
author = {Brezzi, Franco and Marini, Luisa Donatella},
title = {The three-field formulation for elasticity problems},
journal = {GAMM-Mitteilungen},
volume = {28},
number = {2},
pages = {124-153},
doi = {10.1002/gamm.201490016},
year = {2005}
}

@article{NEGRI2015,
title = {Efficient model reduction of parametrized systems by matrix discrete empirical interpolation},
journal = {Journal of Computational Physics},
volume = {303},
pages = {431-454},
year = {2015},
issn = {0021-9991},
doi = {10.1016/j.jcp.2015.09.046},
author = {Federico Negri and Andrea Manzoni and David Amsallem},
}

@article{Joskova_Tyburec_Doskar_2025,
title={Optimising truss-based metamaterials for unimodal and pentamodal behaviour},
volume={54},
DOI={10.14311/APP.2025.54.0034},
 journal={Acta Polytechnica CTU Proceedings}, author={Jošková, Nataša and Tyburec, Marek and Doškář, Martin}, year={2025}, month={Dec.}, pages={34–40} }

@inproceedings{hirschler2026,
  TITLE = {{Optimisation isog{\'e}om{\'e}trique de forme de structures lattices}},
  AUTHOR = {Hirschler, Thibaut and Chasapi, Margarita and Antolin, Pablo and Buffa, Annalisa},
  URL = {https://hal.science/hal-05631873},
  BOOKTITLE = {{17{\`e}me Colloque National en Calcul des Structures}},
  ADDRESS = {Giens, France},
  ORGANIZATION = {{CSMA}},
  HAL_LOCAL_REFERENCE = {113287},
  YEAR = {2026},
  HAL_ID = {hal-05631873},
  HAL_VERSION = {v1},
}

@article{Wang2023,
  author    = {J. Wang and J. Zhu and T. Liu and Y. Wang and H. Zhou and W.-H. Zhang},
  title     = {Topology optimization of gradient lattice structure under harmonic load based on multiscale finite element method},
  journal   = {Structural and Multidisciplinary Optimization},
  year      = {2023},
  volume    = {66},
  number    = {9},
  pages     = {202},
  doi       = {10.1007/s00158-023-03652-3},
  issn      = {1615-147X},
  publisher = {Springer}
}

@article{Zhou2026,
  author    = {H. Zhou and C. Zhou},
  title     = {A multiscale topology optimization design framework with data driven surrogate model},
  journal   = {Scientific Reports},
  year      = {2026},
  volume    = {16},
  pages     = {4647},
  doi       = {10.1038/s41598-025-34920-5},
  publisher = {Springer Nature}
}

@article{KIM2024,
title = {Multiscale topology optimization for the design of spatially-varying three-dimensional lattice structure},
journal = {Computer Methods in Applied Mechanics and Engineering},
volume = {429},
pages = {117140},
year = {2024},
issn = {0045-7825},
doi = {10.1016/j.cma.2024.117140},
author = {Dongjin Kim and Jaewook Lee},
}

\end{document}